\documentclass[]{article}%For paper submission

\usepackage[a4paper,top=2.1cm,bottom=2.60cm,left=2.1cm,right=2.1cm]{geometry} %toglierlo

\usepackage{latexsym}
\usepackage{amsmath}
\usepackage{amsfonts}
\usepackage{mathrsfs}
\usepackage{amstext}
\usepackage{amssymb}
\usepackage{amsthm}       %proof environment

\usepackage{moreverb}

\usepackage[dvips,colorlinks,bookmarksopen,bookmarksnumbered,citecolor=red,urlcolor=red]{hyperref}

\newcommand\BibTeX{{\rmfamily B\kern-.05em \textsc{i\kern-.025em b}\kern-.08em
T\kern-.1667em\lower.7ex\hbox{E}\kern-.125emX}}

\newtheorem{thm}{Theorem}[section]
\newtheorem{prop}[thm]{Proposition}
\newtheorem{defn}[thm]{Definition}

\newtheorem{lem}[thm]{Lemma}

\begin{document}

\title{A singularly perturbed nonlinear traction problem in a periodically perforated domain. A functional analytic approach}

\author{M.~Dalla Riva and P.~Musolino}

\date{}

\maketitle

\noindent
{\bf Abstract:}We consider a periodically perforated  domain obtained by making in $\mathbb{R}^n$ a periodic set of holes, each of them of size proportional to $\epsilon$. Then we introduce a nonlinear boundary value problem for the Lam\'e equations in such a periodically perforated domain. The unknown of the problem is a vector valued function $u$ which represents the displacement attained in the equilibrium configuration by the points of  a periodic linearly elastic matrix with a hole of size $\epsilon$ contained in each periodic cell. We assume that the traction exerted by the matrix on the boundary of each hole depends (nonlinearly) on the displacement attained by the points of the boundary of the hole. Then our aim is to describe what happens to the displacement vector function $u$ when $\epsilon$ tends to $0$. Under suitable assumptions we prove the existence of a family of solutions $\{u(\epsilon,\cdot)\}_{\epsilon\in]0,\epsilon'[}$ with a prescribed limiting behaviour when $\epsilon$ approaches $0$. Moreover, the family $\{u(\epsilon,\cdot)\}_{\epsilon\in]0,\epsilon'[}$ is in a sense locally unique and can be continued real analytically for negative values of $\epsilon$.

\vspace{11pt}

\noindent
{\bf MOS:} 35J65; 31B10; 45F15; 74B05


\vspace{-6pt}

\section{Introduction}\label{introd}
In this article, we consider a singularly perturbed nonlinear traction problem for linearized elastostatics in an infinite periodically perforated domain. We fix once for all
\[
n\in {\mathbb{N}}\setminus\{0,1 \}\,,\qquad  (q_{11},\dots,q_{nn})\in]0,+\infty[^{n}\,.
\]
Here ${\mathbb{N}}$ denotes the 
set of natural numbers including $0$. We denote by $Q$ the fundamental periodicity cell defined by
\begin{equation}\label{Q}
Q\equiv\Pi_{j=1}^{n}]0,q_{jj}[
\end{equation}
 and by $\nu_Q$ the outward unit normal to $\partial Q$, where it exists. We denote by $q$ the diagonal matrix defined by
\begin{equation}\label{q}
q\equiv \left(
\begin{array}{cccc}
q_{11} &   0 & \dots & 0   
\\
0          &q_{22} &\dots & 0
\\
\dots & \dots & \dots & \dots  
\\
0& 0 & \dots & q_{nn}
\end{array}\right)\, .
\end{equation}
%and by $\nu_Q$ the outward unit normal to $\partial Q$, where it exists. 
Then, 
\[
q {\mathbb{Z}}^{n}\equiv  \{qz:\,z\in{\mathbb{Z}}^{n}\}
\]
is the set of vertices of a periodic subdivision of ${\mathbb{R}}^{n}$ corresponding to the fundamental cell $Q$. Let
\[
m\in {\mathbb{N}}\setminus\{0\}\,,\qquad\alpha\in]0,1[\,.
\]
Let $\Omega^h$ be a subset of the Euclidean space $\mathbb{R}^n$ which satisfies the following assumption.
\begin{equation}\label{ass} 
\begin{split}
\text{$\Omega^h$ is a bounded connected open subset of ${\mathbb{R}}^{n}$ of class $C^{m,\alpha}$ such that ${\mathbb{R}}^{n}\setminus{\mathrm{cl}}\Omega^h$ is connected and that $0 \in\Omega^h$}\\
\end{split}
\end{equation}
The letter `h' stands for `hole'. If $p\in Q$ and $\epsilon \in \mathbb{R}$, then we set
\[
\Omega^h_{p,\epsilon} \equiv p+\epsilon\Omega^h \,.
\]
A simple topological argument shows that there exists a real number $\epsilon_0$ such that
\begin{equation}\label{eps0}
\epsilon_0>0\text{ and }\mathrm{cl}\Omega^h_{p,\epsilon} \subseteq Q\text{ for all }\epsilon\in]-\epsilon_0,\epsilon_0[\,.
\end{equation}
Then we denote by  ${\mathbb{S}} [\Omega^h_{p,\epsilon}]^{-}$ the periodically perforated domain defined by
\[
{\mathbb{S}} [\Omega^h_{p,\epsilon}]^{-}\equiv {\mathbb{R}}^{n}\setminus\cup_{z \in \mathbb{Z}^n}{\mathrm{cl}}(\Omega^h_{p,\epsilon}+qz )
\] 
for all $\epsilon\in]-\epsilon_0,\epsilon_0[$.\par

We now introduce a nonlinear traction boundary value problem in ${\mathbb{S}} [\Omega^h_{p,\epsilon}]^{-}$. To do so, we denote by $T$ the function from $ ]1-(2/n),+\infty[\times M_n(\mathbb{R})$ to $M_n(\mathbb{R})$ defined by 
\[
T(\omega,A)\equiv (\omega-1)(\mathrm{tr}A)I_n+(A+A^t) \qquad \forall  \omega \in ]1-(2/n),+\infty[\,,\ A \in M_n(\mathbb{R})\,.
\] Here $M_n(\mathbb{R})$ denotes the space of $n\times n$ matrices with real entries, $I_n$ denotes the $n\times n$ identity matrix, $\mathrm{tr}A$ and $A^t$ denote the trace and the transpose matrix of $A$, respectively. We note that $(\omega-1)$ plays the role of the ratio between the first and second Lam\'e constants and that the classical linearization of the Piola Kirchoff tensor equals the second Lam\'e constant times $T(\omega,\cdot)$ (cf., \textit{e.g.}, Kupradze, Gegelia, Bashele{\u\i}shvili, and Burchuladze \cite{KuGeBaBu79}). Now let $G$ be a (nonlinear) function from $\partial \Omega^h \times\mathbb{R}^n$ to $\mathbb{R}^n$. Let $B \in M_n(\mathbb{R})$. Let $\epsilon \in ]0,\epsilon_0[$. Then we consider the following nonlinear traction boundary value problem 
 \begin{equation}\label{bvp:nltraceleps}
 \left \lbrace 
 \begin{array}{ll}
  \mathrm{div}\, T(\omega, Du )= 0 & \textrm{in ${\mathbb{S}} [\Omega^h_{p,\epsilon}]^-$}\,, \\
u(x+qe_j) =u(x) +Be_j&  \textrm{$\forall x \in \mathrm{cl} {\mathbb{S}} [\Omega^h_{p,\epsilon}]^{-}, \forall j\in \{1,\dots,n\}$}, \\
T(\omega,Du(x))\nu_{\Omega^h_{p,\epsilon}}(x)=G\bigl((x-p)/\epsilon,u(x)\bigr) & \textrm{$\forall x \in \partial \Omega^h_{p,\epsilon}$}\,,
 \end{array}
 \right.
 \end{equation}
where $\{e_1,\dots,e_n\}$ denotes the canonical basis of $\mathbb{R}^n$ and $\nu_{\Omega^h_{p,\epsilon}}$ denotes the outward unit normal to $\partial \Omega^h_{p,\epsilon}$.   We note that due to the presence of a nonlinear term in the third equation we cannot claim the existence of a solution of problem \eqref{bvp:nltraceleps}. However, for a fixed vector $\tilde{\xi} \in \mathbb{R}^n$ and under suitable assumptions we shall prove that there exists $\epsilon'\in]0,\epsilon_0]$ such that problem \eqref{bvp:nltraceleps} has a solution $u(\epsilon,\cdot)$ in $C^{m,\alpha}(\mathrm{cl}\mathbb{S}[{\Omega^h_{p,\epsilon}}]^-,\mathbb{R}^n)$ for all $\epsilon\in]0,\epsilon'[$. The family of solutions $\{u(\epsilon,\cdot)\}_{\epsilon\in]0,\epsilon'[}$ converges to the function $Bq^{-1}(x-p)+\tilde\xi$ of $x\in\mathbb{R}^n$ in a sense which will be clarified in Section \ref{rep}. Moreover, $\{u(\epsilon,\cdot)\}_{\epsilon\in]0,\epsilon'[}$  is unique in a local sense which will be clarified in Section \ref{uniq}.
Then we pose the following questions. 
\begin{enumerate}
\item[(j)] Let $x$ be fixed in $\mathbb{R}^n \setminus (p+q\mathbb{Z}^n)$. What can be said on the map $\epsilon \mapsto u(\epsilon,x)$ when $\epsilon$ is close to $0$ and positive?
\item[(jj)] Let $t$ be fixed in $\mathbb{R}^n \setminus \Omega^h$. What can be said on the map $\epsilon \mapsto u(\epsilon,p+\epsilon t)$ when $\epsilon$ is close to $0$ and positive?
\end{enumerate}
In a sense, question (j) concerns the `macroscopic' behaviour far from the cavities, whereas question (jj) is related to the `microscopic' behaviour of $u(\epsilon,\cdot)$ near the boundary of the holes. 

Questions of this type have long been investigated for linear problems 
with the methods of Asymptotic Analysis and of Calculus of the 
Variations. Thus for example,  one could resort to Asymptotic Analysis and may succeed  to write out an asymptotic  expansion for $u(\epsilon,x)$ and $u(\epsilon,p+\epsilon t)$. In this sense, we mention the work of Ammari and Kang \cite{AmKa07}, Ammari, Kang, and Lee \cite{AmKaLe09}, Ammari, Kang, and Touibi \cite{AmKaTo05}, Ammari, Kang, and Lim \cite{AmKaLi06}, Maz'ya and Movchan \cite{MaMo10}, Maz'ya, Nazarov, and Plamenewskij \cite{MaNaPl00i, MaNaPl00ii}, 
Maz'ya, Movchan, and Nieves \cite{MaMoNi11}. We also mention the extensive literature of Calculus of Variations and of Homogenization Theory,  and in particular the contributions of  
Bakhvalov and Panasenko \cite{BaPa89}, Cioranescu and Murat \cite{CiMu82i, CiMu82ii},
Jikov, Kozlov, and Ole\u{\i}nik \cite{JiKoOl94},
Mar{\v{c}}enko and Khruslov \cite{MaKh74}.\par
 
Furthermore,  boundary value problems in domains with periodic inclusions, for example for the Laplace equation, have been  analysed, at least for the two dimensional case, with the method of functional equations. Here we mention Castro, Pesetskaya, and Rogosin \cite{CaPeRo09},  Drygas and Mityushev \cite{DrMi09}.

In connection with doubly periodic problems for composite materials, we mention the monograph of Grigolyuk  and Fil'shtinskij \cite{GrFi92}.\par

Here we wish to characterize the behaviour of $u(\epsilon,\cdot)$ at $\epsilon=0$ by a different approach. In particular, if we consider a certain function $f(\epsilon)$ relative to the solution $u(\epsilon,\,\cdot\,)$,  as for example one of those in questions (j), (jj) above,  we would try to represent $f(\epsilon)$ for $\epsilon$ small and positive in terms of real analytic maps defined in a whole neighborhood of $\epsilon=0$ and in terms of possibly singular but known functions of $\epsilon$, such as $\epsilon^{-1}$, $\log \epsilon$, etc.. We observe that our approach does have its advantages. Indeed, if for 
example we know that the map in (j) equals for $\epsilon>0$ a real analytic function defined in a whole neighbourhood of $\epsilon=0$, then we know that such a map can be expanded in  power series for $\epsilon$ small. Such a project has been carried out by Lanza de Cristoforis and collaborators in several papers for problems in a bounded domain with one small hole (cf., \textit{e.g.}, \cite{La02, La04,La07a,La08,La10, DaMu12}). For nonlinear problems in the frame of linearized elastostatics, we also mention, \textit{e.g.}, \cite{DaLa10, DaLa10b, DaLa10c, DaLa11}, and for the Stokes equation \cite{Da11}. For problems for the Laplace and Poisson equations in periodically perforated domains, we mention \cite{LaMu11, Mu11, Mu11a, DaMu13}. We note that this paper represents the first step in the analysis of periodic boundary value problems for linearized elastostatics with this approach.

This article is organized as follows. Section \ref{not} is a section of notation and Sections \ref{prel1}, \ref{prel2} are sections of preliminaries. In Section \ref{form} we formulate problem \eqref{bvp:nltraceleps} in terms of an equivalent integral equation which we can analyse by means of the Implicit Function Theorem for real analytic maps. Then we introduce our family of solutions $\{u(\epsilon,\,\cdot\,)\}_{\epsilon\in]0,\epsilon'[}$. In Section \ref{rep}, we prove our main Theorem~\ref{thm:repsol}, where we answer to the questions in (j), (jj). In Section \ref{uniq}, we prove that the family $\{u(\epsilon,\,\cdot\,)\}_{\epsilon\in]0,\epsilon'[}$  is locally unique in a sense which will be clarified.

\section{Some notation}\label{not}

We  denote the norm on 
a   normed space ${\mathcal X}$ by $\|\cdot\|_{{\mathcal X}}$. Let 
${\mathcal X}$ and ${\mathcal Y}$ be normed spaces. We endow the  
space ${\mathcal X}\times {\mathcal Y}$ with the norm defined by 
$\|(x,y)\|_{{\mathcal X}\times {\mathcal Y}}\equiv \|x\|_{{\mathcal X}}+
\|y\|_{{\mathcal Y}}$ for all $(x,y)\in  {\mathcal X}\times {\mathcal 
Y}$, while we use the Euclidean norm for ${\mathbb{R}}^{n}$.
 We denote by $\mathcal{L}(\mathcal{X},\mathcal{Y})$ the space of linear and continuous maps from $\mathcal{X}$ to $\mathcal{Y}$, equipped with its usual norm of the uniform convergence on the unit sphere of $\mathcal{X}$. We denote by $I$ the identity operator. The inverse function of an 
invertible function $f$ is denoted $f^{(-1)}$, as opposed to the 
reciprocal of a real-valued function $g$, or the inverse of a 
matrix $B$, which are denoted $g^{-1}$ and $B^{-1}$, respectively.  For 
standard definitions of Calculus in normed spaces and for the definition and properties of (real) analytic functions in Banach space, we refer to 
Cartan~\cite{Ca71}, Prodi and Ambrosetti~\cite{PrAm73}, Deimling~\cite{De85}. Here we just recall that if $\mathcal{X}$, $\mathcal{Y}$ are (real) Banach spaces, and if $F$ is an operator from an open subset $\mathcal{W}$ of $\mathcal{X}$ to $\mathcal{Y}$, then $F$ is real analytic in $\mathcal{W}$ if for every $x_0 \in \mathcal{W}$ there exist $r>0$ and continuous symmetric $j$-linear operators $A_j$ from $\mathcal{X}^j$ to $\mathcal{Y}$ such that $\sum_{j\geq 1} \|A_j\|_{\mathcal{L}(\mathcal{X}^j,\mathcal{Y})} r^j <\infty$ and $F(x_0+h)=F(x_0)+\sum_{j \geq 1} A_j(h,\dots,h)$ for $\|h\|_{\mathcal{X}} \leq r$ (cf., \textit{e.g.}, Prodi and 
Ambrosetti~\cite[p. 89]{PrAm73} and Deimling \cite[p.~150]{De85}). We note that throughout the paper ``analytic" means ``real analytic".    If $B$ is a 
matrix, then 
 $B_{ij}$ denotes 
the $(i,j)$ entry of $B$. If $x\in\mathbb{R}^n$, then $x_{j}$ denotes the $j$-th coordinate of $x$ and 
$|x|$ denotes the Euclidean modulus of $ x$. A  dot ``$\cdot$'' denotes the inner product in ${\mathbb R}^{n}$. For all $R>0$ and all $x\in{\mathbb{R}}^{n}$ we denote by ${\mathbb{B}}_{n}( x,R)$ the ball $\{
y\in{\mathbb{R}}^{n}:\, | x- y|<R\}$.   If 
$\mathcal{S}$ is a subset of ${\mathbb {R}}^{n}$, then $\mathrm{cl}\mathcal{S}$ 
denotes the 
closure of $\mathcal{S}$ and $\partial\mathcal{S}$ denotes the boundary of $\mathcal{S}$. If we further assume that $\mathcal{S}$ is measurable then $|\mathcal{S}|$ denotes the $n$-dimensional measure of $\mathcal{S}$. Let $q$ be as in definition \eqref{q}. Let $\mathcal{P}$ be a subset of $\mathbb{R}^n$ such that $x+qz \in \mathcal{P}$ for all $x \in \mathcal{P}$ and for all $z \in \mathbb{N}$. We say that a function $f$ on $\mathcal{P}$ is $q$-periodic if
\[
f(x+qz)=f(x) \qquad \forall x \in \mathcal{P}\,,\quad \forall z \in \mathbb{Z}^n\, .
\]
Let $\mathcal{O}$ be an open 
subset of ${\mathbb{R}}^{n}$. Let $k\in\mathbb{N}$.  The space of $k$ times continuously 
differentiable real-valued functions on $\mathcal{O}$ is denoted by 
$C^{k}(\mathcal{O},{\mathbb{R}})$, or more simply by $C^{k}(\mathcal{O})$. If $f\in C^{k}(\mathcal{O})$  then $\nabla f$ denotes the gradient $\left(\frac{\partial f}{\partial
x_1},\dots,\frac{\partial f}{\partial
x_n}\right)$ which we think as a column vector. Let $r\in {\mathbb{N}}\setminus\{0\}$. Let $f\equiv(f_1,\dots,f_r)\in \left(C^{k}(\mathcal{O})\right)^{r}$. Then $Df$ denotes the Jacobian matrix
$\left(\frac{\partial f_s}{\partial
x_l}\right)_{  (s,l)\in\{1,\dots,r\}\times\{1,\dots,n\}}$.  Let  $\eta\equiv
(\eta_{1},\dots ,\eta_{n})\in{\mathbb{N}}^{n}$, $|\eta |\equiv
\eta_{1}+\dots +\eta_{n}  $. Then $D^{\eta} f$ denotes
$\frac{\partial^{|\eta|}f}{\partial
x_{1}^{\eta_{1}}\dots\partial x_{n}^{\eta_{n}}}$.  The
subspace of $C^{k}(\mathcal{O})$ of those functions $f$ whose derivatives $D^{\eta }f$ of
order $|\eta |\leq k$ can be extended with continuity to 
$\mathrm{cl}\mathcal{O}$  is  denoted $C^{k}(
\mathrm{cl}\mathcal{O})$. Let $\beta\in]0,1[$. The
subspace of $C^{k}(\mathrm{cl}\mathcal{O}) $  whose
functions have $k$-th order derivatives that are uniformly
H\"{o}lder continuous in $\mathrm{cl}\mathcal{O}$ with exponent  $\beta$ is denoted $C^{k,\beta} (\mathrm{cl}\mathcal{O})$  
(cf., \textit{e.g.},~Gilbarg and 
Trudinger~\cite{GiTr83}). If $f\in C^{0,\beta}(\mathrm{cl}\mathcal{O})$, then its  $\beta$-H\"{o}lder constant $|f:\mathrm{cl}\mathcal{O}|_{\beta}$ is defined as  
 $\sup\left\{
\frac{|f( x )-f( y)|}{| x- y|^{\beta}
}: x, y\in {\mathrm{cl}}\mathcal{O} ,  x\neq
 y\right\}$. The subspace of $C^{k}(\mathrm{cl}\mathcal{O}) $ of those functions $f$ such that $f_{|{\mathrm{cl}}(\mathcal{O}\cap{\mathbb{B}}_{n}(0,R))}\in
C^{k,\beta}({\mathrm{cl}}(\mathcal{O}\cap{\mathbb{B}}_{n}(0,R)))$ for all $R\in]0,+\infty[$ is denoted $C^{k,\beta}_{{\mathrm{loc}}}(\mathrm{cl}\mathcal{O}) $.  Let 
$\mathcal{S}\subseteq {\mathbb{R}}^{r}$. Then $C^{k
,\beta }(\mathrm{cl}\mathcal{O} ,\mathcal{S})$ denotes
$\left\{f\in \left( C^{k,\beta} (\mathrm{cl}\mathcal{O})\right)^{r}:\ f(
\mathrm{cl}\mathcal{O})\subseteq \mathcal{S}\right\}$. Then we set
\[
\begin{split}
C^{k}_{b}({\mathrm{cl}}\mathcal{O},\mathbb{R}^n)\equiv
\{
u\in C^{k}({\mathrm{cl}}\mathcal{O},\mathbb{R}^n):\,
D^{\eta}u\ {\mathrm{is\ bounded}}\ \textrm{for all }\eta\in {\mathbb{N}}^{n}\
{\mathrm{with}}\ |\eta|\leq k
\}\,,
\end{split}
\]
and we endow $C^{k}_{b}({\mathrm{cl}}\mathcal{O},\mathbb{R}^n)$ with its usual  norm
\[
\|u\|_{ C^{k}_{b}({\mathrm{cl}}\mathcal{O},\mathbb{R}^n) }\equiv
\sum_{\eta \in \mathbb{N}^n\, , \ |\eta|\leq k}\sup_{x\in {\mathrm{cl}}\Omega }|D^{\eta}u(x)|\,.
\]
We define
\[
\begin{split}
C^{k,\beta}_{b}({\mathrm{cl}}\mathcal{O},\mathbb{R}^n)\equiv
\{
u\in C^{k,\beta}({\mathrm{cl}}\mathcal{O},\mathbb{R}^n):\,
D^{\eta}u\ {\mathrm{is\ bounded}}\ \textrm{for all }\eta\in {\mathbb{N}}^{n}\
{\mathrm{with
}}\ |\eta|\leq k
\}\,,
\end{split}
\]
and we endow $C^{k,\beta}_{b}({\mathrm{cl}}\mathcal{O},\mathbb{R}^n)$ with its usual  norm
\[
\|u\|_{ C^{k,\beta}_{b}({\mathrm{cl}}\mathcal{O},\mathbb{R}^n) }\equiv
\sum_{\eta \in \mathbb{N}^n\, ,\ |\eta|\leq k}\sup_{x\in {\mathrm{cl}}\mathcal{O} }|D^{\eta}u(x)|
+\sum_{\eta \in \mathbb{N}^n\, ,\ |\eta| = k}|D^{\eta}u: {\mathrm{cl}}\mathcal{O} |_{\beta}\, .
\]

Let $\mathcal{O} $ be a bounded
open subset of  ${\mathbb{R}}^{n}$. Then $C^{k}(\mathrm{cl}\mathcal{O} )$ 
and $C^{k,\beta}(\mathrm{cl}
\mathcal{O} )$ endowed with their usual norm are well known to be 
Banach spaces  (cf., \textit{e.g.}, Troianiello~\cite[\S 1.2.1]{Tr87}). 
We say that a bounded open subset $\mathcal{O}$ of ${\mathbb{R}}^{n}$ is of class 
$C^{k}$ or of class $C^{k,\beta}$, if its closure is a 
manifold with boundary imbedded in 
${\mathbb{R}}^{n}$ of class $C^{k}$ or $C^{k,\beta}$, respectively
 (cf., \textit{e.g.}, Gilbarg and Trudinger~\cite[\S 6.2]{GiTr83}). 
  For standard properties of functions 
in Schauder spaces, we refer the reader to Gilbarg and 
Trudinger~\cite{GiTr83} and to Troianiello~\cite{Tr87}
(see also Lanza \cite[\S 2, Lem.~3.1, 4.26, Thm.~4.28]{La91}, 
 Lanza  and Rossi \cite[\S 2]{LaRo04}).
If $\mathcal{M}$ is a manifold  imbedded in 
${\mathbb{R}}^{n}$ of class $C^{k,\beta}$ with $k\ge 1$, then one can define the Schauder spaces also on $\mathcal{M}$ by 
exploiting the local parametrization. In particular, if $\mathcal{O}$ is a bounded open set of class $C^{k,\beta}$ with $k\ge 1$, then one can consider 
the space $C^{l,\beta}(\partial\mathcal{O})$ on $\partial\mathcal{O}$
with $l \in \{0,\dots,k\}$ and the trace operator from $C^{l,\beta}({\mathrm{cl}}\mathcal{O})$ to
$C^{l,\beta}(\partial\mathcal{O})$ is linear and continuous.
 Now let $Q$ be as in definition \eqref{Q}. If ${\mathcal{S}_Q}$ is an arbitrary subset of ${\mathbb{R}}^{n}$  such that
 ${\mathrm{cl}}{\mathcal{S}_Q}\subseteq Q$, then we define
\[
{\mathbb{S}} [{\mathcal{S}_Q}]\equiv 
\bigcup_{z\in{\mathbb{Z}}^{n} }(qz+{\mathcal{S}_Q})=q{\mathbb{Z}}^{n}+{\mathcal{S}_Q}\,,\qquad
{\mathbb{S}} [{\mathcal{S}_Q}]^{-}\equiv {\mathbb{R}}^{n}\setminus{\mathrm{cl}}{\mathbb{S}} [{\mathcal{S}_Q}]\,.\nonumber
\] 
We note that if $\mathbb{R}^n\setminus \mathrm{cl}{\mathcal{S}_Q}$ is connected, then $\mathbb{S}[{\mathcal{S}_Q}]^{-}$ is also connected.  If ${\Omega_Q}$ is an open subset of $\mathbb{R}^n$ such that $\mathrm{cl}\Omega_Q \subseteq Q$, then we denote by $C^{k}_{q}({\mathrm{cl}}{\mathbb{S}}[{\Omega_Q}],\mathbb{R}^n )$, $C^{k,\beta}_{q}({\mathrm{cl}}{\mathbb{S}}[{\Omega_Q}],\mathbb{R}^n )$, $C^{k}_{q}({\mathrm{cl}}{\mathbb{S}}[{\Omega_Q}]^-,\mathbb{R}^n )$, and $C^{k,\beta}_{q}({\mathrm{cl}}{\mathbb{S}}[{\Omega_Q}]^-,\mathbb{R}^n )$ the subsets of the $q$-periodic functions belonging to $C^{k}_{b}({\mathrm{cl}}{\mathbb{S}}[{\Omega_Q}],\mathbb{R}^n )$, and to $C^{k,\beta}_{b}({\mathrm{cl}}{\mathbb{S}}[{\Omega_Q}],\mathbb{R}^n )$, and to $C^{k}_{b}({\mathrm{cl}}{\mathbb{S}}[{\Omega_Q}]^-,\mathbb{R}^n )$, and to $C^{k,\beta}_{b}({\mathrm{cl}}{\mathbb{S}}[{\Omega_Q}]^-,\mathbb{R}^n )$, respectively. We regard the sets $C^{k}_{q}({\mathrm{cl}}{\mathbb{S}}[{\Omega_Q}],\mathbb{R}^n )$, $C^{k,\beta}_{q}({\mathrm{cl}}{\mathbb{S}}[{\Omega_Q}],\mathbb{R}^n )$, $C^{k}_{q}({\mathrm{cl}}{\mathbb{S}}[{\Omega_Q}]^-,\mathbb{R}^n )$, and $C^{k,\beta}_{q}({\mathrm{cl}}{\mathbb{S}}[{\Omega_Q}]^-,\mathbb{R}^n )$ as Banach subspaces of $C^{k}_{b}({\mathrm{cl}}{\mathbb{S}}[{\Omega_Q}],\mathbb{R}^n )$, and of $C^{k,\beta}_{b}({\mathrm{cl}}{\mathbb{S}}[{\Omega_Q}],\mathbb{R}^n )$, and of $C^{k}_{b}({\mathrm{cl}}{\mathbb{S}}[{\Omega_Q}]^-,\mathbb{R}^n )$, and of $C^{k,\beta}_{b}({\mathrm{cl}}{\mathbb{S}}[{\Omega_Q}]^-,\mathbb{R}^n )$, respectively. 

\section{Periodic elastic layer potentials}\label{prel1}

We denote by $S_{n}$  the function from  
${\mathbb{R}}^{n}\setminus\{0\}$ to
${\mathbb{R}}$ defined by 
\[
S_{n}(x)\equiv
\left\{
\begin{array}{lll}
\frac{1}{s_{n}}\log |x| \qquad &   \forall x\in 
{\mathbb{R}}^{n}\setminus\{0\},\quad & {\mathrm{if}}\ n=2\,,
\\
\frac{1}{(2-n)s_{n}}|x|^{2-n}\qquad &   \forall x\in 
{\mathbb{R}}^{n}\setminus\{0\},\quad & {\mathrm{if}}\ n>2\,,
\end{array}
\right.
\]
where $s_{n}$ denotes the $(n-1)$-dimensional measure of 
$\partial{\mathbb{B}}_{n}(0,1)$. $S_{n}$ is well-known to be the 
fundamental solution of the Laplace operator.

Let  $\omega \in ]1-(2/n),+\infty[$. We denote by $\Gamma_{n,\omega}(\cdot)$ the matrix valued function from $\mathbb{R}^n \setminus \{0\}$ to $M_{n}(\mathbb{R})$ which takes $x$ to the matrix $\Gamma_{n,\omega}(x)$ with $(i,j)$ entry defined by
\[
\Gamma_{n,\omega,i}^j(x)\equiv \frac{\omega+2}{2(\omega+1)}\delta_{i,j}S_n(x)-\frac{\omega}{2(\omega+1)}\frac{1}{s_n}\frac{x_i x_j}{|x|^n}\qquad\forall (i,j)\in\{1,\dots,n\}^2\,,
\]
where $\delta_{i,j}=1$ if $i=j$, $\delta_{i,j}=0$ if $i \neq j$. As is well known, $\Gamma_{n,\omega}$ is the fundamental solution of the operator
\[
L[\omega]\equiv \Delta+\omega \nabla \mathrm{div}\,.
\]
We note that the classical operator of linearized homogenous isotropic elastostatics equals $L[\omega]$ times the second constant of Lam\'e, and that $L[\omega]u=\mathrm{div} \,T(\omega,Du)$ for all regular vector valued functions $u$, and that the classical fundamental solution of the operator of linearized homogenous and isotropic elastostatics equals $\Gamma_{n,\omega}$ times the reciprocal of the second constant of Lam\'e (cf., \textit{e.g.}, Kupradze, Gegelia, Bashele{\u\i}shvili, and Burchuladze \cite{KuGeBaBu79}). We find also convenient to set
\[
\Gamma_{n,\omega}^j\equiv \bigl(\Gamma_{n,\omega,i}^j\bigr)_{i \in \{1,\dots,n\}}\,,
\]
which we think as a column vector for all $j\in\{1,\dots,n\}$.  

Now let $\alpha \in ]0,1[$. Let $\Omega$ be a bounded open subset of $\mathbb{R}^n$ of class $C^{1,\alpha}$. Then we set
\begin{align}
& v[\omega,\mu](x)\equiv \int_{\partial \Omega}\Gamma_{n,\omega}(x-y)\mu(y)\,d\sigma_y\,,\nonumber\\
& w[\omega,\mu](x)\equiv-\Bigl(\int_{\partial \Omega}\mu^t(y)T(\omega,D \Gamma^i_{n,\omega}(x-y))\nu_{\Omega}(y)\,d\sigma_y\Bigr)_{i\in \{1,\dots,n\}}\,,\nonumber
\end{align}
for all $x \in \mathbb{R}^n$ and for all $\mu \equiv (\mu_j)_{j\in \{1,\dots,n\}} \in C^{0,\alpha}(\partial \Omega,\mathbb{R}^n)$. Here $d\sigma$ denotes the $(n-1)$-dimensional measure on $\partial\Omega$ and $\nu_\Omega$ denotes the outward unit normal to $\partial\Omega$. As is well known, $v[\omega,\mu]$ is continuous in the whole of $\mathbb{R}^n$. We define
\[
v^+[\omega,\mu]\equiv v[\omega,\mu]_{|\mathrm{cl} \Omega}\,, \qquad v^-[\omega,\mu]\equiv v[\omega,\mu]_{|\mathbb{R}^n \setminus \Omega}\,.
\]
Also, $w[\omega,\mu]_{|\Omega}$ admits a unique continuous extension to $\mathrm{cl} \Omega$, which we denote by $w^+[\omega,\mu]$, and $w[\omega,\mu]_{|\mathbb{R}^n \setminus \mathrm{cl}\Omega}$ admits a unique continuous extension to $\mathbb{R}^n \setminus \Omega$, which we denote by $w^-[\omega,\mu]$. We further define
\[
w_{\ast}[\omega, \mu](x)\equiv \int_{\partial \Omega}\sum_{l=1}^n \mu_{l}(y)T(\omega,D\Gamma_{n,\omega}^{l}(x-y))\nu_{\Omega}(x)\,d\sigma_y \qquad \forall x \in \partial \Omega\,,
\]
for all $\mu \equiv (\mu_j)_{j\in \{1,\dots,n\}} \in C^{0,
\alpha}(\partial \Omega,\mathbb{R}^n)$. For properties of elastic layer potentials, we refer, \textit{e.g.}, to \cite[Theorem A.2]{DaLa10}.

In the following Theorem~\ref{psper} we introduce a periodic analogue of the fundamental solution of $L[\omega]$ (cf., \textit{e.g.}, Ammari and Kang \cite[Lemma 9.21]{AmKa07}, Ammari, Kang, and Lim \cite[Lemma 3.2]{AmKaLi06}). To do so we need the following notation.  We denote by ${\mathcal{S}}({\mathbb{R}}^{n},\mathbb{C})$ the Schwartz space of complex valued rapidly decreasing functions. $\mathcal{S}'(\mathbb{R}^n,\mathbb{C})$ denotes the space of complex tempered distributions and $M_{n}\bigl(\mathcal{S}'(\mathbb{R}^n,\mathbb{C})\bigr)$ denotes the set of $n\times n$ matrices with entries in $\mathcal{S}'(\mathbb{R}^n,\mathbb{C})$. The symbols $\bar \zeta$ and $\bar f$ denote the conjugate of a complex number $\zeta$ and of a complex valued function $f$, respectively.  If $y\in{\mathbb{R}}^{n}$ and $f$ is a function defined in ${\mathbb{R}}^{n}$, we set
$\tau_{y}f(x)\equiv f(x-y)$ for all $x\in {\mathbb{R}}^{n}$. If $u\in \mathcal{S}'(\mathbb{R}^n,\mathbb{C})$, then we set
\[
<\tau_{y}u,f>\equiv<u,\tau_{-y}f>\qquad\forall f\in  \mathcal{S}(\mathbb{R}^n,\mathbb{C})\,.
\] Finally, $L^1_{\mathrm{loc}}(\mathbb{R}^n)$ denotes the space of (equivalence classes of) locally summable measurable functions from $\mathbb{R}^n$ to $\mathbb{R}$.

\begin{thm}
\label{psper}
Let $\omega \in ]1-(2/n),+\infty[$. Let $\Gamma_{n,\omega}^{q}\equiv (\Gamma_{n,\omega,j}^{q,k})_{(j,k)\in\{1,\dots,n\}^2}$ be the element of $M_{n}\bigl(\mathcal{S}'(\mathbb{R}^n,\mathbb{C})\bigr)$ with $(j,k)$ entry defined by
\begin{equation}\label{eq:distfs1}
\Gamma_{n,\omega,j}^{q,k}\equiv \sum_{z \in \mathbb{Z}^n \setminus \{0\}} \frac{1}{4 \pi^2 |Q|  |q^{-1}z|^2}\Biggl[ -\delta_{j,k}+\frac{\omega}{\omega+1}\frac{(q^{-1}z)_j(q^{-1}z)_k}{|q^{-1}z|^2}\Biggr]E_{2 \pi iq^{-1}z} \qquad \forall (j,k) \in \{1,\dots,n\}^2\,,
\end{equation}
where $E_{2\pi i q^{-1} z}$ is the function from $\mathbb{R}^n$ to $\mathbb{C}$ defined by
\[
E_{2\pi i q^{-1} z}(x)\equiv e^{2\pi i (q^{-1} z)\cdot x}
\qquad  \forall x\in{\mathbb{R}}^{n}
\] for all $z\in\mathbb{Z}^n$.
 Then the following statements hold.
\begin{enumerate}
\item[(i)]
\[
\tau_{q_{ll}e_l}\Gamma_{n,\omega,j}^{q,k}=\Gamma_{n,\omega,j}^{q,k} \qquad \forall l  \in \{1,\dots,n\}\,,
\]
for all $(j,k) \in \{1,\dots,n\}^2$.
\item[(ii)]
\[
\overline{<\Gamma_{n,\omega,j}^{q,k},\overline{f}>}=<\Gamma_{n,\omega,j}^{q,k},f> \qquad \forall f \in \mathcal{S}(\mathbb{R}^n,\mathbb{C})\,,
\]
for all $(j,k) \in \{1,\dots,n\}^2$.
\item[(iii)]
\[
L[\omega] \Gamma_{n,\omega}^{q}=\sum_{z \in \mathbb{Z}^n}\delta_{qz}I_n-\frac{1}{ |Q|}I_n \qquad \text{$\mathrm{in}$ $M_{n}\bigl(\mathcal{S}'(\mathbb{R}^n,\mathbb{C})\bigr)$}\,,
\]
where $\delta_{qz}$ denotes the Dirac measure with mass at $qz$ for all $z \in \mathbb{Z}^n$.
\item[(iv)] $\Gamma_{n,\omega}^{q}$ is real analytic from $\mathbb{R}^n \setminus q\mathbb{Z}^n$ to $M_n(\mathbb{R})$.
\item[(v)] The difference $\Gamma_{n,\omega}^{q}-\Gamma_{n,\omega}$ can be extended to a real analytic function from $(\mathbb{R}^n \setminus q \mathbb{Z}^n) \cup\{0\}$ to $M_n(\mathbb{R})$ which we denote by $R^q_{n,\omega}$. Moreover
\[
L[\omega] R^q_{n,\omega}=\sum_{z \in \mathbb{Z}^n\setminus\{0\}}\delta_{qz}I_n-\frac{1}{ |Q|}I_n 
\]
in the sense of distributions.
\item[(vi)] $\Gamma_{n,\omega,j}^{q,k} \in L^1_{\mathrm{loc}}(\mathbb{R}^n)$, for all $(j,k)\in \{1,\dots,n\}^2$.
\item[(vii)] $\Gamma_{n,\omega}^q(x)=\Gamma_{n,\omega}^{q}(-x)$ for all $x \in \mathbb{R}^n \setminus q\mathbb{Z}^n$.
\end{enumerate}
\end{thm} 
{\bf Proof.} The Theorem is a simple modification of the corresponding result of \cite[Theorem 3.1]{LaMu10a}, where an analogue of a periodic fundamental solution for a second order strongly elliptic differential operator with constant coefficients has been constructed (see also Ammari and Kang \cite[Lemma 9.21]{AmKa07}, Ammari, Kang, and Lim \cite[Lemma 3.2]{AmKaLi06}). Indeed, since
\[
\sup_{z \in \mathbb{Z}^n \setminus \{0\}}\frac{1}{4 \pi^2 |Q|  |q^{-1}z|^2}\Biggl[-\delta_{j,k}+ \frac{\omega}{\omega+1}\frac{(q^{-1}z)_j(q^{-1}z)_k}{ |q^{-1}z|^2}\Biggr]< +\infty \qquad \forall (j,k) \in \{1,\dots,n\}^2\, ,
\]
 one can prove that the generalized series in \eqref{eq:distfs1} defines a tempered distribution, and accordingly $\Gamma_{n,\omega}^q \in M_{n}(\mathcal{S}'(\mathbb{R}^n,\mathbb{C}))$ (cf.~\cite[Proof of Theorem 3.1]{LaMu10a}). Statement (i) follows by the definition of $\Gamma_{n,\omega}^q$ and by the periodicity of $E_{2\pi i q^{-1} z}$. The statement in (ii) is a straightforward consequence of the obvious equality
\[
\begin{split}
\frac{1}{4 \pi^2 |Q|  |q^{-1}z|^2}\Biggl[ -\delta_{j,k}+\frac{\omega}{\omega+1}\frac{(q^{-1}z)_j(q^{-1}z)_k}{ |q^{-1}z|^2}\Biggr]=\frac{1}{ 4 \pi^2  |Q| |-q^{-1}z|^2}\Biggl[-\delta_{j,k}+ \frac{\omega}{\omega+1}\frac{(-q^{-1}z)_j(-q^{-1}z)_k}{ |-q^{-1}z|^2}\Biggr]\ 
\quad \forall z \in \mathbb{Z}^n\setminus \{0\},
\end{split}
\]
for all $(j,k) \in \{1,\dots,n\}^2$, and of
\[
\overline{<E_{2 \pi i q^{-1}z},\overline{f}>}=<\overline{E_{2 \pi i q^{-1}z}},f>=<E_{2 \pi i q^{-1}(-z)},f> \qquad \forall f \in \mathcal{S}(\mathbb{R}^n,\mathbb{C})\,, \ \forall z \in \mathbb{Z}^n \setminus \{0\}\,. 
\]
We now consider statement (iii). By Poisson's summation formula, we have
\[
\Bigl(\Delta \Gamma_{n,\omega}^q\Bigr)_{jk}=\delta_{j,k}\Biggl[\sum_{z \in \mathbb{Z}^n}\delta_{qz}-\frac{1}{ |Q|}\Biggr]+\sum_{z \in \mathbb{Z}^n \setminus \{0\}}\frac{-1}{ |Q|}\frac{\omega}{\omega+1}\frac{(q^{-1}z)_j(q^{-1}z)_k}{ |q^{-1}z|^2}E_{2 \pi i q^{-1}z}\, ,
\]
and
\[
\begin{split}
\Bigl(\nabla \mathrm{div}\, \Gamma_{n,\omega}^q\Bigr)_{jk}=\sum_{z \in \mathbb{Z}^n \setminus \{0\}}&\Biggl\{\sum_{l=1}^n\Biggl[\frac{1}{4\pi^2 |q^{-1}z|^2 |Q|}\Biggl(-\delta_{l,k}+\frac{\omega}{\omega+1}\frac{(q^{-1}z)_l(q^{-1}z)_k}{ |q^{-1}z|^2}\Biggr)(-4\pi^2)(q^{-1}z)_l(q^{-1}z)_j\Biggr]\Biggr\}E_{2 \pi i q^{-1}z}\\
=\sum_{z \in \mathbb{Z}^n \setminus \{0\}}&\frac{1}{ |q^{-1}z|^2 |Q|}\frac{1}{\omega+1}(q^{-1}z)_j(q^{-1}z)_kE_{2 \pi i q^{-1}z},
\end{split}
\]
for all $(j,k)\in \{1,\dots,n\}^2$. Hence,
\[
\Bigl(\Delta \Gamma_{n,\omega}^q\Bigr)_{jk}+\omega \Bigl(\nabla \mathrm{div}\, \Gamma_{n,\omega}^q\Bigr)_{jk}=\delta_{j,k}\Biggl[\sum_{z \in \mathbb{Z}^n}\delta_{qz}-\frac{1}{ |Q|}\Biggr],
\]
for all $(j,k) \in \{1,\dots,n\}^2$, and thus (iii) follows. Statements (iv), (v) follow by (iii) and by elliptic regularity theory, while (vi) follows by the local integrability of $\Gamma_{n,\omega}$ and the periodicity of $\Gamma_{n,\omega}^q$. Finally, by a straightforward verification based on definition \eqref{eq:distfs1}, statement (vii) easily follows.  Hence, the proof is complete. 
\hfill $\Box$
\vspace{\baselineskip}

We find  convenient to set
\[
\Gamma_{n,\omega}^{q,j}\equiv \bigl(\Gamma_{n,\omega,i}^{q,j}\bigr)_{i \in \{1,\dots,n\}}\,,\qquad
R_{n,\omega}^{q,j}\equiv \bigl(R_{n,\omega,i}^{q,j}\bigr)_{i \in \{1,\dots,n\}}\,,
\]
which we think as column vectors for all $j\in\{1,\dots,n\}$.  

Let $\omega \in ]1-(2/n),+\infty[$. Let  $\alpha\in]0,1[$. Let ${\Omega_Q}$ be a bounded open subset of ${\mathbb{R}}^{n}$ of class $C^{1,\alpha}$ such that ${\mathrm{cl}}{\Omega_Q}\subseteq Q$. Let $\mu \in C^{0,\alpha}(\partial {\Omega_Q},\mathbb{R}^n)$. Then we denote by $v_q[\omega, \mu]$ the periodic single layer potential, namely $v_q[\omega,\mu]$ is the function from $\mathbb{R}^n$ to $\mathbb{R}^n$ defined by 
\[
v_q[\omega, \mu](x)\equiv \int_{\partial {\Omega_Q}}\Gamma^q_{n,\omega}(x-y)\mu(y)\,d\sigma_y \qquad \forall x \in \mathbb{R}^n\,.
\]
We note here that the fundamental solution $\Gamma^q_{n,\omega}$ takes values in $M_n(\mathbb{R})$ (cf.~Theorem \ref{psper} (ii) and (iv)). We also find convenient to set
\begin{equation}\label{eq:vqast}
w_{q,\ast}[\omega, \mu](x)\equiv \int_{\partial {\Omega_Q}}\sum_{l=1}^n \mu_{l}(y)T(\omega,D\Gamma_{n,\omega}^{q,l}(x-y))\nu_{{\Omega_Q}}(x)\,d\sigma_y \quad \forall x \in \partial {\Omega_Q}\,.
\end{equation}

In the following Theorem we collect some properties of the periodic single layer potential.
\begin{thm}
\label{sperpot}
Let  $\alpha\in]0,1[$, $m\in {\mathbb{N}}\setminus\{0\}$. Let ${\Omega_Q}$ be a bounded open subset of ${\mathbb{R}}^{n}$ of class $C^{m,\alpha}$ such that ${\mathrm{cl}}{\Omega_Q}\subseteq Q$. Then the following statements hold.
\begin{enumerate}
\item[(i)] If $\mu\in C^{0,\alpha}(\partial{\Omega_Q},\mathbb{R}^n)$, then $v_{q}[\omega,\mu]$ is $q$-periodic and 
\[
L[\omega]v_{q}[\omega,\mu](x)
=
-\frac{1}{|Q|}\int_{\partial{\Omega_Q}}\mu \,d\sigma
\]
for all  $x\in 
 {\mathbb{R}}^{n}\setminus\partial{\mathbb{S}}[{\Omega_Q}]$.
\item[(ii)]  If $\mu\in C^{m-1,\alpha}(\partial{\Omega_Q},\mathbb{R}^n)$, then the function 
$v^{+}_{q}[\omega,\mu]\equiv v_{q}[\omega,\mu]_{|{\mathrm{cl}}{\mathbb{S}}[{\Omega_Q}]}$ belongs to $C^{m,\alpha}_{q}({\mathrm{cl}}{\mathbb{S}}[{\Omega_Q}],\mathbb{R}^n)$ and the operator
 which takes $\mu$ to 
$v^{+}_{q}[\omega,\mu]  $ is  continuous from $C^{m-1,\alpha}(\partial{\Omega_Q},\mathbb{R}^n)$ to $C^{m,\alpha}_{q}({\mathrm{cl}}
{\mathbb{S}}[{\Omega_Q}],\mathbb{R}^n)$. 
\item[(iii)]  If $\mu\in C^{m-1,\alpha}(\partial{\Omega_Q},\mathbb{R}^n)$, then the function 
$v^{-}_{q}[\omega,\mu]\equiv v_{q}[\omega,\mu]_{|{\mathrm{cl}}{\mathbb{S}}[{\Omega_Q}]^{-}}$ belongs to $C^{m,\alpha}_{q}
({\mathrm{cl}}{\mathbb{S}}[{\Omega_Q}]^{-},\mathbb{R}^n)$ and the operator which takes $\mu$ to $v^{-}_{q}[\omega,\mu]$   is  continuous from  $ C^{m-1,\alpha}(\partial{\Omega_Q},\mathbb{R}^n)$ to $C^{m,\alpha}_{q}
({\mathrm{cl}}{\mathbb{S}}[{\Omega_Q}]^{-},\mathbb{R}^n)$.
\item[(iv)]  The operator which takes $\mu$ to $w_{q,\ast}[\omega,\mu]$ is  continuous from the space $C^{m-1,\alpha}(\partial{\Omega_Q},\mathbb{R}^n)$ to itself, and we have 
\begin{equation}\label{sperpot2}
T\bigl(\omega,Dv_{q}^{\pm}[\omega,\mu](x)\bigr)\nu_{{\Omega_Q}}(x)=\mp\frac{1}{2}\mu(x)+w_{q,\ast}[\omega,\mu](x) \qquad \forall x \in \partial{\Omega_Q}\, ,
\end{equation}
for all $\mu \in C^{m-1,\alpha}(\partial{\Omega_Q},\mathbb{R}^n)$.
\end{enumerate}
\end{thm}
{\bf Proof.} By splitting $\Gamma_{n,\omega}^{q}$ into the sum of $\Gamma_{n,\omega}$ and $R_{n,\omega}^{q}$, by exploiting Theorem \ref{psper} and classical potential theory for linearized elastostatics (cf., \textit{e.g.}, \cite[Theorem A.2]{DaLa10}) and standard properties of integral  operators with real analytic kernels and with no singularity (cf., \textit{e.g.}, \cite[\S 4]{LaMu10b}), one can prove the validity of statements (i), (ii), (iii), and (iv). See also \cite[Theorem 3.7]{LaMu10a}, where the periodic single layer potential for a second order strongly elliptic differential operator with constant coefficients has been constructed. \hfill $\Box$
\vspace{\baselineskip}

Similarly, we introduce the periodic double layer potential $w_q[\omega, \mu]$. 
Let $\omega \in ]1-(2/n),+\infty[$. Let  $\alpha\in]0,1[$. Let ${\Omega_Q}$ be a bounded open subset of ${\mathbb{R}}^{n}$ of class $C^{1,\alpha}$ such that ${\mathrm{cl}}{\Omega_Q}\subseteq Q$. Let $\mu \in C^{0,\alpha}(\partial {\Omega_Q},\mathbb{R}^n)$. We set
\[
w_q[\omega, \mu](x)\equiv -\Biggl(\int_{\partial {\Omega_Q}}\mu^t(y)T(\omega, D \Gamma_{n,\omega}^{q,i}(x-y))\nu_{{\Omega_Q}}(y)\,d\sigma_y\Biggr)_{i\in\{1,\dots,n\}} \  \forall x \in \mathbb{R}^n\,,
\] which we think as a column vector. In the following Theorem we collect some properties of the periodic double layer potential.

\begin{thm}
\label{dperpot}
Let $\omega \in ]1-(2/n),+\infty[$. Let  $\alpha\in]0,1[$, $m\in {\mathbb{N}}\setminus\{0\}$. Let ${\Omega_Q}$ be a bounded  open subset of ${\mathbb{R}}^{n}$ of class $C^{m,\alpha}$ such that ${\mathrm{cl}}{\Omega_Q}\subseteq Q$. Then the following statements hold.
\begin{enumerate}
\item[(i)] If $\mu\in C^{0,\alpha}(\partial{\Omega_Q},\mathbb{R}^n)$, then $w_{q}[\omega,\mu]$ is $q$-periodic and
\[
L[\omega]w_{q}[\omega,\mu](x)=0\qquad\forall x\in {\mathbb{R}}^{n}\setminus\partial{\mathbb{S}}[{\Omega_Q}]\,.
\]
\item[(ii)] If $\mu\in C^{m,\alpha}(\partial{\Omega_Q},\mathbb{R}^n)$, then 
the restriction $w_{q}[\omega,\mu]_{|{\mathbb{S}}[{\Omega_Q}]}$ can be extended to a function $w^{+}_{ q }[\omega,\mu] \in C_{q}^{m,\alpha}({\mathrm{cl}}{\mathbb{S}}[{\Omega_Q}],\mathbb{R}^n)$, and the  restriction $w_{q}[\omega,\mu]_{|{\mathbb{S}}[{\Omega_Q}]^{-}}$ can be extended to a function $w_{ q }^{-}[\omega,\mu] \in C^{m,\alpha}_{q}({\mathrm{cl}}{\mathbb{S}}[{\Omega_Q}]^{-},\mathbb{R}^n)$, and we have
\begin{equation}\label{dperpot2a}
w^{\pm }_{q}[\omega,\mu]=\pm\frac{1}{2}\mu+
w_{q}[\omega,\mu]
\qquad{\mathrm{on}}\ \partial {\Omega_Q} \,.
\end{equation}
\item[(iii)] The operator from  $C^{m,\alpha}(\partial{\Omega_Q},\mathbb{R}^n)$ to 
$C^{m,\alpha}_{q}({\mathrm{cl}}{\mathbb{S}}[{\Omega_Q}],\mathbb{R}^n)$ which takes $\mu$ to $
w_{q}^{+}[\omega,\mu]$ is  continuous. The operator  from   $C^{m,\alpha}(\partial{\Omega_Q},\mathbb{R}^n)$ to $C^{m,\alpha}_{q}({\mathrm{cl}}{\mathbb{S}}[{\Omega_Q}]^{-},\mathbb{R}^n)$ which takes $\mu$ to $w_{q}^{-}[\omega,\mu]$ is continuous. 
\item[(iv)] We have
\begin{equation}\label{dbperpot3}
w_{q}[\omega,e_j](x)=
\left\{
\begin{array}{lll}
-\frac{|{\Omega_Q}|}{|Q|}e_j  & \text{if $x \in \mathbb{S}[{\Omega_Q}]^{-}$} \, ,\\
\biggl(\frac{1}{2}-\frac{|{\Omega_Q}|}{|Q|}\biggr) e_j &   \text{if $x \in \partial \mathbb{S}[{\Omega_Q}]$} \, ,
\\
\biggl(1-\frac{|{\Omega_Q}|}{|Q|}\biggr)e_j& \text{if $x \in \mathbb{S}[{\Omega_Q}]$} \, ,
\end{array}
\right.
\end{equation}
for all $j \in \{1,\dots,n\}$.
\end{enumerate}
\end{thm}
{\bf Proof.} 
By splitting $\Gamma_{n,\omega}^{q}$ into the sum of $\Gamma_{n,\omega}$ and $R_{n,\omega}^{q}$, by exploiting Theorem \ref{psper} and classical potential theory for linearized elastostatics (cf., \textit{e.g.}, \cite[Theorem A.2]{DaLa10}) and standard properties of integral  operators with real analytic kernels and with no singularity (cf., \textit{e.g.}, \cite[\S 4]{LaMu10b}), one can prove the validity of statements (i), (ii), and (iii). See also \cite[Theorem 3.18]{LaMu10a}, where the periodic double layer potential for a second order strongly elliptic differential operator with constant coefficients has been constructed. We now turn to the proof of statement (iv). It clearly suffices to prove equality \eqref{dbperpot3} for $x \in \mathbb{S}[{\Omega_Q}]^{-}$. Indeed, case $x \in \partial \mathbb{S}[{\Omega_Q}]$ and case $x \in \mathbb{S}[{\Omega_Q}]$ can be proved by exploiting the case $x \in \mathbb{S}[{\Omega_Q}]^{-}$ and the jump relations of equality \eqref{dperpot2a}. By periodicity, we can assume $x \in \mathrm{cl}Q \setminus \mathrm{cl}{\Omega_Q}$. By the Divergence Theorem and Theorem \ref{psper} (iii), we have 
\[
\begin{split}
-\int_{\partial {\Omega_Q}}e_j^t T(\omega, D \Gamma_{n,\omega}^{q,i}(x-y))\nu_{{\Omega_Q}}(y)\,d\sigma_y=\int_{{\Omega_Q}}e_j^t\Bigl(L[\omega] \Gamma_{n,\omega}^{q,i}(x-y)\Bigr) \,dy =-\frac{|{\Omega_Q}|}{|Q|} \delta_{i,j}\,,
\end{split}
\]
for all $(i,j) \in \{1,\dots,n\}^2$. As a consequence, statement (iv) follows. Thus the proof is complete.
\hfill $\Box$
\vspace{\baselineskip}

\section{Some preliminary results on periodic problems for linearized elastostatics}\label{prel2}

In the following Propositions~\ref{prop:uneu} and \ref{prop:ossint} we consider a periodic boundary value problem for linearized elastostatics and we show some properties of the corresponding solution.

\begin{prop}\label{prop:uneu}
Let $\omega \in ]1-(2/n),+\infty[$. Let  $\alpha\in]0,1[$, $m\in {\mathbb{N}}\setminus\{0\}$. Let ${\Omega_Q}$ be a bounded open subset of ${\mathbb{R}}^{n}$ of class $C^{m,\alpha}$ such that 
${\mathbb{R}}^{n}\setminus{\mathrm{cl}}{\Omega_Q}$ is connected and that ${\mathrm{cl}}{\Omega_Q}\subseteq Q$. Let $u \in  C^{m,\alpha}_q(\mathrm{cl} \mathbb{S}[{\Omega_Q}]^-,\mathbb{R}^n)$ be a solution of
\[
\left \lbrace 
 \begin{array}{ll}
L[\omega]u(x)= 0 & \textrm{$\forall x \in \mathbb{S}[{\Omega_Q}]^-$}\,, \\
u(x+qe_k) =u(x) &  \textrm{$\forall x \in \mathrm{cl} \mathbb{S}[{\Omega_Q}]^-\,, \quad \forall k \in \{1,\dots,n\}$}\,, \\
T(\omega,Du(x))\nu_{{\Omega_Q}}(x)=0 & \textrm{$\forall x \in \partial {\Omega_Q}$}\,.
 \end{array}
 \right.
\]
Then there exists $b \in \mathbb{R}^n$ such that $u(x)=b$ for all $x \in \mathrm{cl} \mathbb{S}[{\Omega_Q}]^-$.
\end{prop}
{\bf Proof.} By the periodicity of $u$ we have $\int_{\partial Q}u^tT(\omega,Du)\nu_Q\,d\sigma=0$. Thus the Divergence Theorem implies that
\[
\int_{ Q\setminus \mathrm{cl}{\Omega_Q}}\mathrm{tr} \bigl(T(\omega,Du)D^tu\bigr)\,dx=-\int_{\partial {\Omega_Q}}u^tT(\omega,Du)\nu_{\Omega_Q} \,d\sigma=0\,.
\]
Then $\mathrm{tr} \bigl(T(\omega,Du)D^tu\bigr)=0$ in $Q\setminus \mathrm{cl}{\Omega_Q}$, and by arguing as in \cite[Proposition 2.1]{DaLa10}, one can prove that there exist a skew symmetric matrix $A \in M_{n}(\mathbb{R})$ and $b \in \mathbb{R}^n$, such that
\[
u(x)=Ax+b \qquad \forall x \in \mathrm{cl}  Q\setminus \mathrm{cl}{\Omega_Q}\,.
\]
By the periodicity of $u$, we have
\[
A qe_k=u(qe_k)-u(0)=0 \qquad \forall k \in \{1,\dots,n\}\,.
\]
Accordingly, $A=0$. Hence, $u(x)=b$ for all $x \in \mathrm{cl}  Q\setminus \mathrm{cl}{\Omega_Q}$, and thus, by periodicity, $u(x)=b$ for all $x \in \mathrm{cl} \mathbb{S}[{\Omega_Q}]^-$.
\hfill $\Box$
\vspace{\baselineskip}

\begin{prop}\label{prop:ossint}
Let $\omega \in ]1-(2/n),+\infty[$. Let  $\alpha\in]0,1[$, $m\in {\mathbb{N}}\setminus\{0\}$. Let ${\Omega_Q}$ be a bounded open subset of ${\mathbb{R}}^{n}$ of class $C^{m,\alpha}$ such that ${\mathrm{cl}}{\Omega_Q}\subseteq Q$. Let $u \in  C^{m,\alpha}_q(\mathrm{cl} \mathbb{S}[{\Omega_Q}]^-,\mathbb{R}^n)$ be such that
\[
\left \lbrace 
 \begin{array}{ll}
L[\omega]u(x)= 0 & \textrm{$\forall x \in \mathbb{S}[{\Omega_Q}]^-$}\,, \\
u(x+qe_k) =u(x) &  \textrm{$\forall x \in \mathrm{cl} \mathbb{S}[{\Omega_Q}]^-, \quad \forall k \in \{1,\dots,n\}$}\,. \\
 \end{array}
 \right.
\]
Then 
\[
\int_{\partial {\Omega_Q}}T(\omega,Du)\nu_{{\Omega_Q}}\,d\sigma=0\,.
\]
\end{prop}
{\bf Proof.}  By the periodicity of $u$ we have $\int_{\partial Q}T(\omega,Du)\nu_Q\,d\sigma=0$. Then, by the Divergence Theorem one verifies that 
\[
\int_{\partial {\Omega_Q}} T(\omega, Du(y))\nu_{{\Omega_Q}}(y)\,d\sigma_y=-\int_{ Q\setminus \mathrm{cl}{\Omega_Q}}\mathrm{div}\, \Bigl(T\bigl(\omega, Du(y)\bigr)\Bigr)\,dy=0\, , 
\]
and the conclusion follows.
\hfill $\Box$
\vspace{\baselineskip}

In Proposition~\ref{prop:bij} below, we show that $\frac{1}{2}I+w_{q,\ast}[\omega,\cdot]$ is a homeomorphism from $C^{m-1,\alpha}(\partial\Omega_Q,\mathbb{R}^n)$ to itself (cf.~definition \eqref{eq:vqast}). To do so we need the following technical Lemma~\ref{lmm:intdens}.

\begin{lem} \label{lmm:intdens}
Let $\omega \in ]1-(2/n),+\infty[$. Let  $\alpha\in]0,1[$, $m\in {\mathbb{N}}\setminus\{0\}$. Let ${\Omega_Q}$ be a bounded open subset of ${\mathbb{R}}^{n}$ of class $C^{m,\alpha}$ such that ${\mathrm{cl}}{\Omega_Q}\subseteq Q$. Let $\mu \in C^{0,\alpha}(\partial {\Omega_Q},\mathbb{R}^n)$. Then 
\[
\int_{\partial {\Omega_Q}}w_{q,\ast}[\omega,\mu]\,d\sigma= \left(\frac{1}{2}-\frac{|{\Omega_Q}|}{|Q|}\right)\int_{\partial {\Omega_Q}}\mu\,d\sigma\,.
\]
\end{lem}
{\bf Proof.} By the properties of the composition of ordinary and singular integrals, and by Theorems \ref{psper} (vii) and \ref{dperpot} (iv) we have
\[
\begin{split}
\int_{\partial {\Omega_Q}}\int_{\partial {\Omega_Q}}\sum_{l=1}^n \mu_{l}(y)T(\omega,D\Gamma_{n,\omega}^{q,l}(x-y))\nu_{{\Omega_Q}}(x)\,d\sigma_y \,d\sigma_x&= \int_{\partial {\Omega_Q}}\sum_{l=1}^n \mu_{l}(y)\int_{\partial {\Omega_Q}}T(\omega,D\Gamma_{n,\omega}^{q,l}(x-y))\nu_{{\Omega_Q}}(x)\,d\sigma_x \,d\sigma_y\\
&
= \int_{\partial {\Omega_Q}}\sum_{l=1}^n \mu_{l}(y)\Biggl(\frac{1}{2}-\frac{|{\Omega_Q}|}{|Q|}\Biggr)e_l\,d\sigma_y  = \left(\frac{1}{2}-\frac{|{\Omega_Q}|}{|Q|}\right)\int_{\partial {\Omega_Q}}\mu\,d\sigma\,,
\end{split}
\] 
and thus the proof is complete.
\hfill $\Box$
\vspace{\baselineskip}

Then we have the following.

\begin{prop}\label{prop:bij}
Let $\omega \in ]1-(2/n),+\infty[$. Let  $\alpha\in]0,1[$, $m\in {\mathbb{N}}\setminus\{0\}$. Let ${\Omega_Q}$ be a bounded open subset of ${\mathbb{R}}^{n}$ of class $C^{m,\alpha}$ such that 
${\mathbb{R}}^{n}\setminus{\mathrm{cl}}{\Omega_Q}$ is connected and that ${\mathrm{cl}}{\Omega_Q}\subseteq Q$.Then $\frac{1}{2}I +w_{q,\ast}[\omega,\cdot]$ is a linear homeomorphism from $C^{m-1,\alpha}(\partial {\Omega_Q},\mathbb{R}^n)$ to itself.
\end{prop}
{\bf Proof.} We observe that $\frac{1}{2}I +w_{q,\ast}[\omega,\cdot]$ is a continuous linear operator from $C^{m-1,\alpha}(\partial {\Omega_Q},\mathbb{R}^n)$ to itself (cf.~Theorem \ref{sperpot} (iv)). Thus by the Open Mapping Theorem, in order to prove that $\frac{1}{2}I +w_{q,\ast}[\omega,\cdot]$ is an homeomorphism, it suffices to show that it is a bijection. To do so, we verify that $\frac{1}{2}I +w_{q,\ast}[\omega,\cdot]$ is a Fredholm operator of index $0$ from $C^{m-1,\alpha}(\partial {\Omega_Q},\mathbb{R}^n)$ to itself and has null space $\{0\}$. Let $\mu \in C^{m-1,\alpha}(\partial {\Omega_Q},\mathbb{R}^n)$. We have
\[
w_{q,\ast}[\omega,\mu](x)=w_{\ast}[\omega,\mu](x)+\int_{\partial {\Omega_Q}}\sum_{l=1}^n \mu_{l}(y)T(\omega,D R_{n,\omega}^{q,l}(x-y))\nu_{{\Omega_Q}}(x)\,d\sigma_y \qquad  \forall x \in \partial {\Omega_Q}\,.
\]
Since $R_{n,\omega,i}^{q,j}(\cdot)$ is real analytic in $(\mathbb{R}^n\setminus q\mathbb{Z}^n)\cup\{0\}$ for all $(i,j)\in \{1,\dots,n\}^2$, standard properties of integral  operators with real analytic kernels and with no singularity (cf., \textit{e.g.}, \cite{LaMu10b}), the compactness of the embedding of $C^{m,\alpha}(\partial {\Omega_Q},\mathbb{R}^n)$ into $C^{m-1,\alpha}(\partial {\Omega_Q},\mathbb{R}^n)$, and standard calculus in Schauder spaces imply that the map from $C^{m-1,\alpha}(\partial {\Omega_Q},\mathbb{R}^n)$ to itself, which takes $\mu$ to the function from $\partial {\Omega_Q}$ to $\mathbb{R}^n$, defined by
\[
\int_{\partial {\Omega_Q}}\sum_{l=1}^n \mu_{l}(y)T(\omega,D R_{n,\omega}^{q,l}(x-y))\nu_{{\Omega_Q}}(x)\,d\sigma_y \qquad \forall x \in \partial {\Omega_Q}\,,
\]
is compact. Since $\frac{1}{2}I +w_{\ast}[\omega,\cdot]$ is a Fredholm operator of index $0$ from $C^{m-1,\alpha}(\partial {\Omega_Q},\mathbb{R}^n)$ to itself (cf., \textit{e.g.}, \cite[Theorem A.9]{DaLa10}), and since compact perturbations of Fredholm operators of index $0$ are Fredholm operators of index $0$, we conclude that $\frac{1}{2}I +w_{q,\ast}[\omega,\cdot]$ is a Fredholm operator of index $0$ from $C^{m-1,\alpha}(\partial {\Omega_Q},\mathbb{R}^n)$ to itself. Now  let $\mu \in C^{m-1,\alpha}(\partial {\Omega_Q},\mathbb{R}^n)$ be such that
\begin{equation}\label{eq:bij1}
\frac{1}{2}\mu(x)+w_{q,\ast}[\omega,\mu](x)=0 \qquad \forall x \in \partial {\Omega_Q}\,.
\end{equation}
We have
\[
T(\omega,Dv_q^-[\omega,\mu])\nu_{{\Omega_Q}}=0 \qquad \text{on $\partial {\Omega_Q}$}
\]
(cf.~equality~\eqref{sperpot2}). Moreover, by Lemma \ref{lmm:intdens} and by equality \eqref{eq:bij1}, we deduce that $\int_{\partial {\Omega_Q}}\mu \,d\sigma=0$. By Theorem \ref{sperpot} and Proposition \ref{prop:uneu}, there exists $b \in \mathbb{R}^n$ such that $v_q^-[\omega,\mu]=b$ in $\mathrm{cl} \mathbb{S}[{\Omega_Q}]^-$. Since
\[
v_q^+[\omega,\mu]=v_q^-[\omega,\mu]=b \qquad \text{on $\partial {\Omega_Q}$}\,,
\]
and by uniqueness results for the Dirichlet problem for $L[\omega]$ in ${\Omega_Q}$, we have $v_q^+[\omega,\mu]=b$ in $\mathrm{cl} {\Omega_Q}$. As a consequence,
\[
T(\omega,Dv_q^+[\omega,\mu])\nu_{{\Omega_Q}}=0 \qquad \text{on $\partial {\Omega_Q}$}\,.
\]
Thus,
\[
\mu = T(\omega,Dv_q^-[\omega,\mu])\nu_{{\Omega_Q}}-T(\omega,Dv_q^+[\omega,\mu])\nu_{{\Omega_Q}} =0 \qquad \text{on $\partial {\Omega_Q}$}
\] 
(cf.~equality~\eqref{sperpot2}). Hence, $\frac{1}{2}I +w_{q,\ast}[\omega,\cdot]$ is an injective Fredholm operator of index $0$ and accordingly a linear homeomorphism from $C^{m-1,\alpha}(\partial {\Omega_Q},\mathbb{R}^n)$ to itself.
\hfill $\Box$
\vspace{\baselineskip}

In the following Proposition \ref{prop:rep} we show a representation formula for a periodic function $u$ defined on the set $\mathrm{cl}\mathbb{S}[{\Omega_Q}]^{-}$ and such that $L[\omega]u=0$. To do so we need the following notation. If   ${\Omega}$ is a bounded open subset of ${\mathbb{R}}^{n}$ of class $C^{m,\alpha}$, with  $\alpha\in]0,1[$, $m\in {\mathbb{N}}\setminus\{0\}$, then we set
\[
C^{m-1,\alpha}(\partial {\Omega},\mathbb{R}^n)_0\equiv \left\{f \in C^{m-1,\alpha}(\partial {\Omega},\mathbb{R}^n)\colon \int_{\partial {\Omega}}f \, d\sigma=0\right\}\,.
\]

\begin{prop}\label{prop:rep}
Let $\omega \in ]1-(2/n),+\infty[$. Let  $\alpha\in]0,1[$, $m\in {\mathbb{N}}\setminus\{0\}$. Let ${\Omega_Q}$ be a bounded open subset of ${\mathbb{R}}^{n}$ of class $C^{m,\alpha}$ such that 
${\mathbb{R}}^{n}\setminus{\mathrm{cl}}{\Omega_Q}$ is connected and such that ${\mathrm{cl}}{\Omega_Q}\subseteq Q$. Let $u \in C^{m,\alpha}_q(\mathrm{cl}\mathbb{S}[{\Omega_Q}]^-,\mathbb{R}^n)$. Assume that
\[
L[\omega]u(x)=0 \qquad \forall x \in \mathbb{S}[{\Omega_Q}]^-\,.
\]
Then there exists a unique pair $(\mu,b) \in C^{m-1,\alpha}(\partial {\Omega_Q},\mathbb{R}^n)_0\times \mathbb{R}^n$ such that
\begin{equation}\label{eq:rep1}
u(x)=v^-_q[\omega,\mu](x)+b \qquad \forall x \in \mathrm{cl}\mathbb{S}[{\Omega_Q}]^-\,.
\end{equation}
\end{prop}
{\bf Proof.} By Proposition~\ref{prop:bij} there exists a unique function $\mu\in C^{m-1,\alpha}(\partial {\Omega_Q},\mathbb{R}^n)$ such that 
\begin{equation}\label{eq:rep2}
\frac{1}{2}\mu(x)+w_{q,\ast}[\omega,\mu](x)=T(\omega,Du(x))\nu_{{\Omega_Q}}(x)  \qquad  \forall x \in \partial {\Omega_Q}\,.
\end{equation} Then Proposition~\ref{prop:ossint} and Lemma~\ref{lmm:intdens} imply that $\int_{\partial{\Omega_Q}}\mu\, d\sigma=0$. Thus $\mu$ belongs to $C^{m-1,\alpha}(\partial {\Omega_Q},\mathbb{R}^n)_0$. By Theorem~\ref{sperpot}, by equation \eqref{eq:rep2}, and by Proposition~\ref{prop:uneu} there exists a unique $b\in\mathbb{R}^n$ such that 
equality \eqref{eq:rep1} holds. Hence there exists a unique pair $(\mu,b)$ in $C^{m-1,\alpha}(\partial {\Omega_Q},\mathbb{R}^n)_0\times \mathbb{R}^n$  such that equality \eqref{eq:rep1} holds.
\hfill $\Box$
\vspace{\baselineskip}

\section{Formulation of an auxiliary problem in terms of an integral equation}\label{form}

In this Section, we convert problem \eqref{bvp:nltraceleps} in the unknown $u$, into an equivalent auxiliary problem. Then we shall provide a formulation of the auxiliary problem in terms of an integral equation.

To do so, we introduce the following notation. Let $m \in \mathbb{N} \setminus \{0\}$, $\alpha \in ]0,1[$. Let $\Omega^h$ be as in assumption \eqref{ass}. If $G\in C^0(\partial \Omega^h \times \mathbb{R}^n,\mathbb{R}^n)$, then we denote by $F_G$ the (nonlinear nonautonomous) composition operator from $C^0(\partial \Omega^h,\mathbb{R}^n)$ to itself which takes $v \in C^0(\partial \Omega^h,\mathbb{R}^n)$ to the function $F_G[v]$ from $\partial \Omega^h$ to $\mathbb{R}^n$, defined by
\[
F_G[v](t)\equiv G(t,v(t)) \qquad \forall t \in \partial \Omega^h\, .
\]
Then we consider the following assumptions.
\begin{align}
& G \in C^0(\partial \Omega^h \times \mathbb{R}^n,\mathbb{R}^n)\,.\label{assG}\\
&\text{$F_G$ maps $C^{m-1,\alpha}(\partial \Omega^h,\mathbb{R}^n)$ to itself.}\label{assFG}
\end{align}
We also note here that if $G \in C^0(\partial \Omega^h\times\mathbb{R}^n,\mathbb{R}^n)$ is such that $F_G$ is real analytic from $C^{m-1,\alpha}(\partial \Omega^h,\mathbb{R}^n)$ to itself, then one can prove that the gradient matrix $D_uG(\cdot,\cdot)$ of $G(\cdot,\cdot)$ with respect to the variable in $\mathbb{R}^n$ exists. If $\tilde{v} \in C^{m-1,\alpha}(\partial \Omega^h,\mathbb{R}^n)$ and $dF_G[\tilde{v}]$ denotes the Fr\'{e}chet differential of $F_G$ at $\tilde{v}$, then we have
\begin{equation}\label{dFG}
dF_{G}[\tilde{v}](v)=
\sum_{l=1}^{n}F_{ \partial_{u_{l}} G  }
[\tilde{v}] v_{l}
\qquad
\forall v\in 
 C^{m-1,\alpha}(\partial\Omega^{h},{\mathbb{R}}^{n})
\end{equation}
(cf.~Lanza \cite[Prop. 6.3]{La07a}). Moreover, 
\begin{equation}\label{DuG}
D_uG(\cdot,\xi)\in C^{m-1,\alpha}(\partial \Omega^h,M_{n}(\mathbb{R})) \qquad \forall \xi \in \mathbb{R}^n\,,
\end{equation}
where $C^{m-1,\alpha}(\partial \Omega^h,M_{n}(\mathbb{R}))$ denotes the space of functions of class $C^{m-1,\alpha}$ from $\partial \Omega^h$ to $M_n(\mathbb{R})$. 

Now let $p\in Q$. Let $\epsilon_0$ be as in assumption \eqref{eps0}. Let $B\in M_n(\mathbb{R})$. Let assumption \eqref{assG} hold. Let $\epsilon \in ]0,\epsilon_0[$. Then one verifies that a function $u \in C^{m,\alpha}_{\mathrm{loc}}(\mathrm{cl}\mathbb{S}[{\Omega^h_{p,\epsilon}}]^-,\mathbb{R}^n)$ solves problem \eqref{bvp:nltraceleps}, if and only if the function $u_\#$ defined by
\[
u_\#(x)\equiv u(x)-Bq^{-1}x \quad \forall x \in \mathrm{cl}\mathbb{S}[{\Omega^h_{p,\epsilon}}]^-\,
\]
is a solution of the following auxiliary problem
 \begin{equation}\label{bvp:auxnltraceleps}
 \left \lbrace 
 \begin{array}{ll}
L[\omega]u_\#(x)= 0 &\forall x \in \mathbb{S} [{\Omega^h_{p,\epsilon}}]^-\,, \\
u_\#(x+qe_j) =u_\#(x)&  \forall x \in \mathrm{cl} {\mathbb{S}} [{\Omega^h_{p,\epsilon}}]^{-},\  \forall  j\in \{1,\dots,n\}\,, \\
T(\omega,Du_\#(x))\nu_{{\Omega^h_{p,\epsilon}}}(x)+T(\omega,Bq^{-1})\nu_{{\Omega^h_{p,\epsilon}}}(x)=G\bigl((x-p)/\epsilon,u_\#(x)+Bq^{-1}x\bigr) & \textrm{$\forall x \in \partial {\Omega^h_{p,\epsilon}}$}\,.
 \end{array}
 \right.
 \end{equation}

We shall now transform the auxiliary problem \eqref{bvp:auxnltraceleps} into an integral equation, by exploiting the representation formula of Proposition \ref{prop:rep} with ${\Omega_Q}$ replaced by ${\Omega^h_{p,\epsilon}}$. We  note that the representation formula of Proposition \ref{prop:rep} includes integrations on the $\epsilon$-dependent domain $\partial {\Omega^h_{p,\epsilon}}$. In order to get rid of such a dependence, we  introduce the following Lemma~\ref{lem:rep}, where we properly rescale the density of the representation formula of Proposition \ref{prop:rep}.

\begin{lem}\label{lem:rep}
Let $\omega \in ]1-(2/n),+\infty[$. Let  $\alpha\in]0,1[$, $m\in {\mathbb{N}}\setminus\{0\}$. Let $\Omega^h$ be as in assumption \eqref{ass}.  Let $p\in Q$. Let $\epsilon_0$ be as in assumption \eqref{eps0}.  Let $\epsilon \in ]0,\epsilon_0[$. Let $u \in C^{m,\alpha}_q(\mathrm{cl}\mathbb{S}[{\Omega^h_{p,\epsilon}}]^-,\mathbb{R}^n)$ be such that
\[
L[\omega]u(x)=0 \qquad \forall x \in \mathbb{S}[{\Omega^h_{p,\epsilon}}]^-\,.
\]
Then there exists a unique pair $(\theta,b) \in C^{m-1,\alpha}(\partial \Omega^h,\mathbb{R}^n)_0\times \mathbb{R}^n$ such that
\[
u(x)=\epsilon^{n-1}\int_{\partial \Omega^h}\Gamma_{n,\omega}^q(x-p-\epsilon s)\theta(s)\,d\sigma_s+b \qquad \forall x \in \mathrm{cl}\mathbb{S}[{\Omega^h_{p,\epsilon}}]^-\,.
\]
\end{lem}
{\bf Proof.} It is a straightforward consequence of Proposition \ref{prop:rep}, of the Theorem of change of variables in integrals, and of standard properties of functions in Schauder spaces.
\hfill $\Box$
\vspace{\baselineskip}

We are now ready to transform problem \eqref{bvp:auxnltraceleps} into an integral equation by means of the following.

\begin{prop}\label{prop:biju}
Let $\omega \in ]1-(2/n),+\infty[$. Let  $\alpha\in]0,1[$, $m\in {\mathbb{N}}\setminus\{0\}$. Let $\Omega^h$ be as in assumption \eqref{ass}.   Let $p\in Q$. Let $\epsilon_0$ be as in assumption \eqref{eps0}. Let $B\in M_n(\mathbb{R})$. Let $G$ be as in assumptions \eqref{assG}, \eqref{assFG}.  Let $\Lambda$ be the map from $]-\epsilon_0,\epsilon_0[\times C^{m-1,\alpha}(\partial \Omega^h,\mathbb{R}^n)_0\times \mathbb{R}^n$ to $C^{m-1,\alpha}(\partial \Omega^h,\mathbb{R}^n)$, defined by
\[
\begin{split}
\Lambda[\epsilon,\theta,\xi](t)
\equiv& \frac{1}{2}\theta(t)+w_{\ast}[\omega,\theta](t)+\epsilon^{n-1} \int_{\partial \Omega^h}\sum_{l=1}^n \theta_{l}(s)T(\omega,D R_{n,\omega}^{q,l}(\epsilon(t-s)))\nu_{\Omega^h}(t)\,d\sigma_s +T(\omega,{Bq^{-1}})\nu_{\Omega^h}(t)\\
&-G\Bigl(t,\epsilon v[\omega,\theta](t)+\epsilon^{n-1}\int_{\partial \Omega^h}R_{n,\omega}^q(\epsilon(t-s))\theta(s)\,d\sigma_s+\epsilon Bq^{-1}t + \xi\Bigr) \qquad  \qquad \qquad \forall t \in \partial \Omega^h\,,
\end{split}
\]
for all $(\epsilon,\theta,\xi)\in ]-\epsilon_0,\epsilon_0[\times C^{m-1,\alpha}(\partial \Omega^h,\mathbb{R}^n)_0\times \mathbb{R}^n$.
If $\epsilon \in ]0,\epsilon_0[$, then the map $u_\#[\epsilon,\cdot,\cdot]$ from the set of pairs $(\theta,\xi)\in C^{m-1,\alpha}(\partial \Omega^h,\mathbb{R}^n)_0 \times \mathbb{R}^n$ that solve the equation
\begin{equation}
\begin{split}\label{eq:biju1}
\Lambda[\epsilon,\theta,\xi]=0
\end{split}
\end{equation}
to the set of functions $u_\# \in C^{m,\alpha}_q(\mathrm{cl} \mathbb{S}[{\Omega^h_{p,\epsilon}}]^-,\mathbb{R}^n)$ which solve problem \eqref{bvp:auxnltraceleps}, which takes $(\theta,\xi)$ to 
\begin{equation}\label{eq:biju2}
u_\#[\epsilon,\theta,\xi](x)\equiv \epsilon^{n-1}\int_{\partial \Omega^h}\Gamma_{n,\omega}^q(x-p-\epsilon s)\theta(s)\,d\sigma_s-Bq^{-1}p+\xi\quad\forall x\in \mathrm{cl}\mathbb{S}[\Omega^h_ {p,\epsilon}]^-\,,
\end{equation}
is a bijection.
\end{prop}
{\bf Proof.} Let $\epsilon \in]0,\epsilon_0[$. Assume that the function $u_\#$ in $C^{m,\alpha}_q(\mathrm{cl} \mathbb{S}[{\Omega^h_{p,\epsilon}}]^-,\mathbb{R}^n)$ solves problem \eqref{bvp:auxnltraceleps}. Then by Lemma \ref{lem:rep}, there exists a unique pair $(\theta,\xi)$ in $C^{m-1,\alpha}(\partial \Omega^h,\mathbb{R}^n)_0\times \mathbb{R}^n$ such that $u_\#$ equals the right hand side of definition \eqref{eq:biju2}. Then a simple computation based on the Theorem of change of variables in integrals and on Theorem \ref{sperpot}, shows that the pair $(\theta,\xi)$ must solve equation \eqref{eq:biju1}. Conversely, one can easily show that if the pair $(\theta,\xi)$ of $C^{m-1,\alpha}(\partial \Omega^h,\mathbb{R}^n)_0\times \mathbb{R}^n$ solves equation \eqref{eq:biju1}, then the function delivered by definition \eqref{eq:biju2} is a solution of problem \eqref{bvp:auxnltraceleps}.
\hfill $\Box$
\vspace{\baselineskip}

Hence we are reduced to analyse equation \eqref{eq:biju1}. We note that for $\epsilon =0$ we obtain an equation which we address to as the \textit{limiting equation} and which has the following form
\begin{equation}\label{eq:lim1}
\begin{split}
\frac{1}{2}\theta(t)+w_{\ast}[\omega,\theta](t)+T(\omega,{Bq^{-1}})\nu_{\Omega^h}(t)-G(t,\xi)=0 \qquad \forall t \in \partial \Omega^h\,.
\end{split}
\end{equation}
Then we have the following Proposition, which shows, under suitable assumptions, the solvability of the limiting equation.

\begin{prop}\label{prop:limsys}
Let $\omega \in ]1-(2/n),+\infty[$. Let  $\alpha\in]0,1[$, $m\in {\mathbb{N}}\setminus\{0\}$. Let $\Omega^h$ be as in assumption \eqref{ass}.    Let $B\in M_n(\mathbb{R})$. Let $G$ be as in assumptions \eqref{assG}, \eqref{assFG}. Assume that there exists $\tilde{\xi} \in \mathbb{R}^n$ such that
\[
\int_{\partial \Omega^h}G(t,\tilde{\xi})\, d\sigma_t=0\,.
\]
Then the integral equation
\[
\frac{1}{2}\theta(t)+w_{\ast}[\omega,\theta](t)+T(\omega,{Bq^{-1}})\nu_{\Omega^h}(t)-G(t,\tilde{\xi})=0 \qquad \forall t \in \partial \Omega^h
\]
has a unique solution in $C^{m-1,\alpha}(\partial \Omega^h,\mathbb{R}^n)_0$, which we denote by $\tilde{\theta}$. As a consequence, the pair $(\tilde{\theta},\tilde{\xi})$ is a solution in $C^{m-1,\alpha}(\partial \Omega^h,\mathbb{R}^n)_0\times \mathbb{R}^n$ of the limiting equation \eqref{eq:lim1}.
\end{prop}
{\bf Proof.} A simple computation based on the Divergence Theorem shows that
\[
\int_{\partial \Omega^h}T(\omega,{Bq^{-1}})\nu_{\Omega^h}(t)\, d\sigma_t =0\, .
\]
Then the Proposition follows by \cite[Remark A.8 and equality (A.7)]{DaLa10}. \hfill $\Box$
\vspace{\baselineskip}

In Theorem \ref{thm:Lmbd} below, we analyse equation \eqref{eq:biju1} around the degenerate value $\epsilon=0$. To do so, we need the following result of classical potential theory for linearized elastostatics.

\begin{prop}\label{prop:bijC}
Let  $\alpha\in]0,1[$, $m\in {\mathbb{N}}\setminus\{0\}$. Let $\Omega$ be a bounded open subset of ${\mathbb{R}}^{n}$ of class $C^{m,\alpha}$ such that 
${\mathbb{R}}^{n}\setminus{\mathrm{cl}}\Omega$ is connected. Let $C\equiv (c_{ij}(\cdot))_{(i,j) \in \{1,\dots, n\}^2} \in C^{m-1,\alpha}(\partial \Omega,M_{n}(\mathbb{R}))$ be such that the matrix
\[
\int_{\partial \Omega}C(y)\,d\sigma_y \equiv \Biggl(\int_{\partial \Omega}c_{ij}(y)\,d\sigma_y\Biggr)_{(i,j)\in \{1,\dots,n\}^2}
\]
is invertible. Then the map from $ C^{m-1,\alpha}(\partial \Omega,\mathbb{R}^n)_0\times \mathbb{R}^n$ to $C^{m-1,\alpha}(\partial \Omega,\mathbb{R}^n)$, which takes $(\mu,b)$ to the function
\[
\frac{1}{2}\mu+w_{\ast}[\omega,\mu]+C b\,,
\] 
is a linear homeomorphism.
\end{prop}
{\bf Proof.} Let $\mathcal{H}$ be the map from $ C^{m-1,\alpha}(\partial \Omega,\mathbb{R}^n)_0\times \mathbb{R}^n$ to $C^{m-1,\alpha}(\partial \Omega,\mathbb{R}^n)$ defined by
\[
\mathcal{H}[\mu,b](x)\equiv\frac{1}{2}\mu(x)+w_{\ast}[\omega,\mu](x)+C(x) b \qquad \forall x \in \partial \Omega\,,
\]
for all $(\mu,b) \in  C^{m-1,\alpha}(\partial \Omega,\mathbb{R}^n)_0\times \mathbb{R}^n$. By standard properties of elastic layer potentials, $\mathcal{H}$ is linear and continuous (cf., \textit{e.g.}, \cite[Theorem A.2]{DaLa10}). Thus, by the Open Mapping Theorem, it suffices to prove that it is a bijection. So let $\psi \in C^{m-1,\alpha}(\partial \Omega,\mathbb{R}^n)$. We need to prove that there exists a unique pair $(\mu,b) \in  C^{m-1,\alpha}(\partial \Omega,\mathbb{R}^n)_0\times \mathbb{R}^n$ such that
\begin{equation}\label{eq:bijC1}
\frac{1}{2}\mu(x)+w_{\ast}[\omega,\mu](x)+C(x) b=\psi(x) \qquad \forall x \in \partial \Omega\,.
\end{equation}
We first prove uniqueness. Let us assume that the pair $(\mu,b) \in  C^{m-1,\alpha}(\partial \Omega,\mathbb{R}^n)_0\times \mathbb{R}^n$ solve equation \eqref{eq:bijC1}. By integrating both sides of equation \eqref{eq:bijC1}, and by the well known identity
\begin{equation}\label{eq:bijC0}
\int_{\partial \Omega}\Bigl(\frac{1}{2}\mu(x)+w_{\ast}[\omega,\mu](x)\Bigr)\,d\sigma_x=\int_{\partial \Omega}\mu(x)\,d\sigma_x
\end{equation}
(cf., \textit{e.g.}, \cite[equality (A.7)]{DaLa10}), we obtain
\[
\Bigl(\int_{\partial \Omega}C(x)\,d\sigma_x\Bigr)b=\int_{\partial \Omega}\psi(x)\,d\sigma_x\,,
\]
and thus
\begin{equation}\label{eq:bijC2}
b=\Bigl(\int_{\partial \Omega}C(x)\,d\sigma_x\Bigr)^{-1}\int_{\partial \Omega}\psi(x)\,d\sigma_x\,.
\end{equation}
As a consequence, $\mu$ is the unique solution in $C^{m-1,\alpha}(\partial \Omega,\mathbb{R}^n)$ of equation
\begin{equation}\label{eq:bijC3}
\frac{1}{2}\mu(x)+w_{\ast}[\omega,\mu](x)=\psi(x)-C(x) \Bigl(\int_{\partial \Omega}C(y)\,d\sigma_y\Bigr)^{-1}\int_{\partial \Omega}\psi(y)\,d\sigma_y \qquad \forall x \in \partial \Omega
\end{equation}
(cf.~\cite[Remark A.8]{DaLa10}). We also note that by equality \eqref{eq:bijC0} the unique solution of equation \eqref{eq:bijC3} is in $ C^{m-1,\alpha}(\partial \Omega,\mathbb{R}^n)_0$. Hence uniqueness follows. In order to prove existence, it suffices to observe that the pair $(\mu,b) \in  C^{m-1,\alpha}(\partial \Omega,\mathbb{R}^n)_0\times \mathbb{R}^n$ identified by equations \eqref{eq:bijC2}, \eqref{eq:bijC3} solves equation \eqref{eq:bijC1} (see also \cite[Remark A.8]{DaLa10}).
\hfill $\Box$
\vspace{\baselineskip}

\begin{thm}\label{thm:Lmbd}
Let $\omega \in ]1-(2/n),+\infty[$. Let  $\alpha\in]0,1[$, $m\in {\mathbb{N}}\setminus\{0\}$. Let $\Omega^h$ be as in assumption \eqref{ass}.   Let $p\in Q$. Let $\epsilon_0$ be as in assumption \eqref{eps0}. Let $B\in M_n(\mathbb{R})$. Let $G$ be as in assumption \eqref{assG}.  Assume that
\begin{equation}\label{assFGan}
\text{$F_G$ is real analytic from $C^{m-1,\alpha}(\partial \Omega^h,\mathbb{R}^n)$ to itself.}
\end{equation}
Assume that there exists $\tilde{\xi} \in \mathbb{R}^n$ such that
\begin{equation}\label{assxi}
\text{$\int_{\partial \Omega^h}G(t,\tilde{\xi})\, d\sigma_t=0$\quad and \quad $\mathrm{det}\int_{\partial \Omega^h}D_uG(t,\tilde{\xi})\, d\sigma_t\neq0$.}
\end{equation}
Let $\Lambda$ be as in Proposition \ref{prop:biju}. Let $\tilde{\theta}$ be the unique function in $C^{m-1,\alpha}(\partial \Omega^h,\mathbb{R}^n)_0$ such that $\Lambda[0,\tilde\theta,\tilde{\xi}]=0$ (cf.~Proposition \ref{prop:limsys}). Then there exist $\epsilon_1 \in ]0,\epsilon_0]$, an open neighbourhood $\mathcal{U}$ of $(\tilde{\theta},\tilde{\xi})$ in $C^{m-1,\alpha}(\partial \Omega^h,\mathbb{R}^n)_0\times\mathbb{R}^n$, and a real analytic map $(\Theta,\Xi)$ from $]-\epsilon_1,\epsilon_1[$ to $\mathcal{U}$, such that the set of zeros of the map $\Lambda$ in $]-\epsilon_1,\epsilon_1[\times\mathcal{U}$ coincides with the graph of $(\Theta,\Xi)$. In particular, $(\Theta[0],\Xi[0])=(\tilde{\theta},\tilde{\xi})$.
\end{thm}
{\bf Proof.} We plan to apply the Implicit Function Theorem for real analytic maps. We first prove that $\Lambda$ is real analytic from $]-\epsilon_0,\epsilon_0[\times C^{m-1,\alpha}(\partial \Omega^h,\mathbb{R}^n)_0\times\mathbb{R}^n$ to $C^{m-1,\alpha}(\partial \Omega^h,\mathbb{R}^n)$. 
We note that
\[
\epsilon \partial\Omega^h-\epsilon \partial \Omega^h \subseteq (\mathbb{R}^n \setminus q\mathbb{Z}^n)\cup \{0\} \qquad \forall \epsilon \in ]-\epsilon_0,\epsilon_0[\,.
\]
Then by standard properties of integral  operators with real analytic kernels and with no singularity (cf., \textit{e.g.}, \cite[\S 4]{LaMu10b}) we can deduce the analyticity of the map from $]-\epsilon_0,\epsilon_0[\times C^{m-1,\alpha}(\partial \Omega^h,\mathbb{R}^n)_0$ to $C^{m-1,\alpha}(\partial \Omega^h,\mathbb{R}^n)$ which takes $(\epsilon,\theta)$ to the function
\[
 \epsilon^{n-1} \int_{\partial \Omega^h}\sum_{l=1}^n \theta_{l}(s)T(\omega,D R_{n,\omega}^{q,l}(\epsilon(t-s)))\nu_{\Omega^h}(t)\,d\sigma_s
\]
of the variable $t \in \partial \Omega^h$. Similarly, the map from $]-\epsilon_0,\epsilon_0[\times C^{m-1,\alpha}(\partial \Omega^h,\mathbb{R}^n)_0$ to $C^{m-1,\alpha}(\partial \Omega^h,\mathbb{R}^n)$ which takes $(\epsilon,\theta)$ to the function
\[
 \epsilon^{n-1}\int_{\partial \Omega^h}R_{n,\omega}^q(\epsilon(t-s))\theta(s)\,d\sigma_s\,
\]
of the variable $t \in \partial \Omega^h$ is real analytic. By classical potential theory for linearized elastostatics, $v[\omega,\cdot]$ and $w_{\ast}[\omega,\cdot]$ are linear and continuous maps from $C^{m-1,\alpha}(\partial \Omega^h,\mathbb{R}^n)$ to $C^{m,\alpha}(\partial\Omega^h,\mathbb{R}^n)$ and to $C^{m-1,\alpha}(\partial \Omega^h,\mathbb{R}^n)$, respectively (cf., \textit{e.g.}, \cite[Theorem A.2]{DaLa10}). Then by standard calculus in Banach spaces and assumption \eqref{assFGan}, we deduce that $\Lambda$ is real analytic from $]-\epsilon_0,\epsilon_0[\times C^{m-1,\alpha}(\partial \Omega^h,\mathbb{R}^n)_0\times\mathbb{R}^n$ to $C^{m-1,\alpha}(\partial \Omega^h,\mathbb{R}^n)$. By standard calculus in Banach space, the differential of $\Lambda$ at $(0,\tilde{\theta},\tilde{\xi})$ with respect to the variables $(\theta,\xi)$ is delivered by the following formula
\[
\begin{split}
\partial_{(\theta,\xi)}\Lambda[0,\tilde{\theta},\tilde{\xi}](\theta^\sharp,\xi^\sharp)(t)=\frac{1}{2}\theta^\sharp(t)+w_{\ast}[\omega,\theta^\sharp](t)-D_uG(t,\tilde{\xi})\xi^\sharp \quad \forall t \in \partial \Omega^h\,,
\end{split} 
\]
for all $(\theta^\sharp,\xi^\sharp)\in C^{m-1,\alpha}(\partial \Omega^h,\mathbb{R}^n)_0\times \mathbb{R}^n$ (see also formula \eqref{dFG}). By assumption \eqref{assxi} and Proposition \ref{prop:bijC}, we deduce that $\partial_{(\theta,\xi)}\Lambda[0,\tilde{\theta},\tilde{\xi}]$ is a linear homeomorphism from $C^{m-1,\alpha}(\partial \Omega^h,\mathbb{R}^n)_0\times \mathbb{R}^n$ to $C^{m-1,\alpha}(\partial \Omega^h,\mathbb{R}^n)$ (see also \eqref{DuG}). Then in order to conclude the proof it suffices to apply the Implicit Function Theorem for real analytic maps in Banach spaces (cf., \textit{e.g.}, Prodi and Ambrosetti \cite[Theorem 11.6]{PrAm73}, Deimling \cite[Theorem 15.3]{De85}).
\hfill $\Box$
\vspace{\baselineskip}

We are now in the position to introduce the following.

\begin{defn}\label{def:ueps}
Let the notation and assumptions of Theorem \ref{thm:Lmbd} hold. Let $u_\#[\cdot,\cdot,\cdot]$ be as in Proposition \ref{prop:biju}. Then we set
\[
u(\epsilon,x)\equiv u_\#[\epsilon,\Theta[\epsilon],\Xi[\epsilon]](x)+Bq^{-1}x \qquad \forall x \in \mathrm{cl} \mathbb{S}[{\Omega^h_{p,\epsilon}}]^-\,,\ \forall \epsilon \in ]0,\epsilon_1[.
\]
\end{defn}

We note that for each $\epsilon \in ]0,\epsilon_1[$ the function $u(\epsilon,\cdot)$ of Definition~\ref{def:ueps} is a solution of problem \eqref{bvp:nltraceleps}.

\section{A functional analytic representation theorem for the family $\{u(\epsilon,\cdot)\}_{\epsilon \in ]0,\epsilon_1[}$}\label{rep}

In the following Theorem \ref{thm:repsol} we show that the family of functions $\{u(\epsilon,\cdot)\}_{\epsilon \in ]0,\epsilon_1[}$ introduced in Definition~\ref{def:ueps} can be continued real analytically for negative values of $\epsilon$, and we answer to the questions in (j), (jj) of the Introduction. 

\begin{thm}\label{thm:repsol}
Let $\omega \in ]1-(2/n),+\infty[$. Let  $\alpha\in]0,1[$, $m\in {\mathbb{N}}\setminus\{0\}$. Let $\Omega^h$ be as in assumption \eqref{ass}.   Let $p\in Q$. Let $\epsilon_0$ be as in assumption \eqref{eps0}. Let $B\in M_n(\mathbb{R})$. Let $G$ be as in assumptions \eqref{assG}, \eqref{assFGan}. Let $\tilde\xi\in\mathbb{R}^n$. Let assumption \eqref{assxi} hold. Then the following statements hold. 
\begin{enumerate}
\item[(i)]Let $\tilde{\Omega}$ be a bounded open subset of $\mathbb{R}^n$ such that $\mathrm{cl} \tilde{\Omega} \subseteq \mathbb{R}^n \setminus (p+q\mathbb{Z}^n)$. Let $k \in \mathbb{N}$. Then there exist $\tilde{\epsilon} \in ]0,\epsilon_0]$ and a real analytic operator $U$ from $]-\tilde{\epsilon},\tilde{\epsilon}[$ to $C^k(\mathrm{cl}\tilde{\Omega},\mathbb{R}^n)$, such that 
\begin{equation} \label{eq:repsol1a0}
\mathrm{cl} \tilde{\Omega} \subseteq \mathbb{S}[{\Omega^h_{p,\epsilon}}]^- \qquad \forall \epsilon \in ]-\tilde{\epsilon},\tilde{\epsilon}[\, ,
\end{equation}
and that
\begin{equation}\label{eq:repsol1a}
u(\epsilon,x)= U[\epsilon](x)\qquad \forall x \in \mathrm{cl} \tilde{\Omega}\, , \quad \forall\epsilon \in]0,\tilde{\epsilon}[\,.
\end{equation}
Moreover,
\begin{equation}\label{eq:repsol1b}
U[0](x)=Bq^{-1}(x-p)+\tilde{\xi}\qquad \forall x \in \mathrm{cl}\tilde{\Omega}\,.
\end{equation}
\item[(ii)]Let $\tilde{\Omega}_r$ be a bounded open subset of $\mathbb{R}^n \setminus \mathrm{cl} \Omega^h$. Then there exist $\tilde{\epsilon}_r \in ]0,\epsilon_0]$ and a real analytic operator $U_{r}$ from $]-\tilde{\epsilon}_r,\tilde{\epsilon}_r[$ to $C^{m,\alpha}(\mathrm{cl}\tilde{\Omega}_r,\mathbb{R}^n)$, such that
\[
p + \epsilon \mathrm{cl} \tilde{\Omega}_r \subseteq Q \setminus {\Omega^h_{p,\epsilon}}\qquad \forall \epsilon \in ]-\tilde{\epsilon}_r,\tilde{\epsilon}_r[\setminus \{0\}\, ,
\]
and that
\begin{equation}\label{eq:repsol2a}
u(\epsilon,p+\epsilon t)=U_{r}[\epsilon](t)\qquad \forall t \in \mathrm{cl} \tilde{\Omega}_r\, ,\quad \forall \epsilon \in]0,\tilde{\epsilon}_r[\,.
\end{equation}
Moreover,
\begin{equation}\label{eq:repsol2b}
U_{r}[0](t)=\tilde{\xi}\qquad \forall t \in \mathrm{cl}\tilde{\Omega}_r\,.
\end{equation}
(Here the letter `r' stands for `rescaled'.)
\end{enumerate}
\end{thm}
 {\bf Proof.} Let $\epsilon_1$, $\Theta$, $\Xi$ be as in Theorem \ref{thm:Lmbd}. We start by proving (i). By taking $\tilde{\epsilon} \in ]0,\epsilon_1]$ small enough, we can assume that condition \eqref{eq:repsol1a0} holds. Consider now equality \eqref{eq:repsol1a}. If $\epsilon \in ]0,\tilde{\epsilon}[$, a simple computation based on the Theorem of change of variables in integrals shows that
 \[
 u(\epsilon,x)=\epsilon^{n-1}\int_{\partial \Omega^h}\Gamma_{n,\omega}^q(x-p-\epsilon s)\Theta[\epsilon](s)\,d\sigma_s-{Bq^{-1}}p+\Xi[\epsilon]+{Bq^{-1}}x \qquad \forall x \in \mathrm{cl}\tilde{\Omega}\,.
 \]
 Thus it is natural to set
 \[
U[\epsilon](x)\equiv\epsilon^{n-1}\int_{\partial \Omega^h}\Gamma_{n,\omega}^q(x-p-\epsilon s)\Theta[\epsilon](s)\,d\sigma_s-{Bq^{-1}}p+\Xi[\epsilon]+{Bq^{-1}}x\qquad \forall x \in \mathrm{cl}\tilde{\Omega}\, ,
 \]
 for all $\epsilon \in ]-\tilde{\epsilon},\tilde{\epsilon}[$. Then we note that
\[
\mathrm{cl}\tilde{\Omega}-p-\epsilon \partial \Omega^h \subseteq \mathbb{R}^n \setminus q\mathbb{Z}^n \qquad \forall \epsilon \in ]-\tilde{\epsilon},\tilde{\epsilon}[\,.
\]
As a consequence, by standard properties of integral  operators with real analytic kernels and with no singularity (cf., \textit{e.g.}, \cite[\S 3]{LaMu10b}), we can conclude that the map from $]-\tilde{\epsilon},\tilde{\epsilon}[$ to $C^k(\mathrm{cl}\tilde{\Omega},\mathbb{R}^n)$, which takes $\epsilon$ to the function
\[
\epsilon^{n-1}\int_{\partial \Omega^h}\Gamma_{n,\omega}^q(x-p-\epsilon s)\Theta[\epsilon](s)\,d\sigma_s
\]
of the variable $x \in \mathrm{cl}\tilde{\Omega}$, is real analytic. Accordingly, $U$ is real analytic from $]-\tilde{\epsilon},\tilde{\epsilon}[$ to $C^{k}(\mathrm{cl}\tilde{\Omega},\mathbb{R}^n)$. By the definition of $U$, equality \eqref{eq:repsol1a} holds. Moreover, the validity of equality \eqref{eq:repsol1b} is obvious, and so the proof of (i) is complete. 

We now consider (ii). Let $R>0$ be such that $(\mathrm{cl}\tilde{\Omega}_r\cup \mathrm{cl} \Omega^h)\subseteq \mathbb{B}_n(0,R)$. By the continuity of the restriction operator from $C^{m,\alpha}(\mathrm{cl}\mathbb{B}_n(0,R)\setminus \Omega^h,\mathbb{R}^n)$ to $C^{m,\alpha}(\mathrm{cl}\tilde{\Omega}_r,\mathbb{R}^n)$, it suffices to prove statement (ii) with $\tilde{\Omega}_r$ replaced by $\mathbb{B}_n(0,R)\setminus \mathrm{cl}\Omega^h$. By taking $\tilde{\epsilon}_r \in ]0,\epsilon_1]$ small enough, we can assume that 
\[
p+\epsilon \mathrm{cl}\mathbb{B}_n(0,R)\subseteq Q \qquad \forall \epsilon \in ]-\tilde{\epsilon}_r,\tilde{\epsilon}_r[\,.
\] 
If $\epsilon \in ]0,\tilde{\epsilon}_r[$, a simple computation based on the Theorem of change of variables in integrals shows that
 \[
 \begin{split}
 u(\epsilon,p+\epsilon t)=&\epsilon \int_{\partial \Omega^h}\Gamma_{n,\omega}(t-s)\Theta[\epsilon](s)\,d\sigma_s+\epsilon^{n-1}\int_{\partial \Omega^h}R_{n,\omega}^q(\epsilon(t- s))\Theta[\epsilon](s)\,d\sigma_s\\
  &-{Bq^{-1}}p+\Xi[\epsilon]+{Bq^{-1}}p+\epsilon {Bq^{-1}}t\qquad  \qquad \qquad \qquad \forall t \in \mathrm{cl}\mathbb{B}_n(0,R)\setminus \Omega^h\,.
 \end{split}
 \]
 Thus it is natural to set
 \[
U_{r}[\epsilon](t)\equiv\epsilon\int_{\partial \Omega^h}\Gamma_{n,\omega}(t-s)\Theta[\epsilon](s)\,d\sigma_s+\epsilon^{n-1}\int_{\partial \Omega^h}R_{n,\omega}^q(\epsilon(t- s))\Theta[\epsilon](s)\,d\sigma_s+\Xi[\epsilon]+\epsilon {Bq^{-1}}t \qquad \forall t \in \mathrm{cl}\mathbb{B}_n(0,R)\setminus \Omega^h\, ,
 \]
 for all $\epsilon \in ]-\tilde{\epsilon}_r,\tilde{\epsilon}_r[$. We note that
 \[
 U_r[\epsilon](t)=\epsilon v^-[\omega,\Theta[\epsilon]](t)+\tilde{U}_{r}[\epsilon](t)+ \Xi[\epsilon]+\epsilon {Bq^{-1}}t \qquad \forall t \in \mathrm{cl}\mathbb{B}_n(0,R)\setminus \Omega^h\, ,
 \]
 for all $\epsilon \in ]-\tilde{\epsilon}_r,\tilde{\epsilon}_r[$, where 
\[
\tilde{U}_{r}[\epsilon](t)\equiv \epsilon^{n-1}\int_{\partial\Omega^h}R_{n,\omega}^{q}(\epsilon(t-s))\Theta[\epsilon](s)\,d\sigma_{s} \qquad \forall t \in \mathrm{cl}\mathbb{B}_n(0,R) \,,
\]
for all $\epsilon \in ]-\tilde{\epsilon}_r,\tilde{\epsilon}_r[$. Then we observe that
\[
\epsilon\mathrm{cl}\mathbb{B}_n(0,R)-\epsilon \partial \Omega^h \subseteq (\mathbb{R}^n \setminus q\mathbb{Z}^n)\cup \{0\} \qquad \forall \epsilon \in ]-\tilde{\epsilon}_r,\tilde{\epsilon}_r[\,.
\]
Accordingly, by standard properties of integral  operators with real analytic kernels and with no singularity (cf., \textit{e.g.}, \cite[\S 4]{LaMu10b}), we can conclude that $\tilde{U}_{r}$ is real analytic from $]-\tilde{\epsilon}_r,\tilde{\epsilon}_r[$ to $C^{m,\alpha}(\mathrm{cl}\mathbb{B}_n(0,R),\mathbb{R}^n)$. By classical results of potential theory and by the real analyticity of $\Theta$, the map from $]-\tilde{\epsilon}_r,\tilde{\epsilon}_r[$ to $C^{m,\alpha}(\mathrm{cl}\mathbb{B}_n(0,R) \setminus \Omega^h,\mathbb{R}^n)$, which takes $\epsilon$ to $v^-[\omega,\Theta[\epsilon]]_{|\mathrm{cl}\mathbb{B}_n(0,R) \setminus \Omega^h}$ is real analytic (cf., \textit{e.g.}, \cite[Theorem A.2]{DaLa10}). Then by the continuity of the restriction operator from $C^{m,\alpha}(\mathrm{cl}\mathbb{B}_n(0,R),\mathbb{R}^n)$ to $C^{m,\alpha}(\mathrm{cl}\mathbb{B}_n(0,R)\setminus \Omega^h,\mathbb{R}^n)$, we deduce that $U_{r}$ is a real analytic map from $]-\tilde{\epsilon}_r,\tilde{\epsilon}_r[$ to $C^{m,\alpha}(\mathrm{cl}\mathbb{B}_n(0,R)\setminus \Omega^h,\mathbb{R}^n)$ and satisfies equalities \eqref{eq:repsol2a}, \eqref{eq:repsol2b} with $\tilde{\Omega}_r$ replaced by $\mathbb{B}_n(0,R)\setminus \mathrm{cl}\Omega^h$. 
 \hfill $\Box$
\vspace{\baselineskip}

\section{Local uniqueness of the family $\{u(\epsilon,\cdot)\}_{\epsilon \in ]0,\epsilon_1[}$}\label{uniq}

In this Section, we show that the family $\{u(\epsilon,\cdot)\}_{\epsilon \in ]0,\epsilon_1[}$ is essentially unique. Namely, we have the following.
\begin{thm}\label{thm:uniqueps}
Let $\omega \in ]1-(2/n),+\infty[$. Let  $\alpha\in]0,1[$, $m\in {\mathbb{N}}\setminus\{0\}$. Let $\Omega^h$ be as in assumption \eqref{ass}.   Let $p\in Q$. Let $\epsilon_0$ be as in assumption \eqref{eps0}. Let $B\in M_n(\mathbb{R})$. Let $G$ be as in assumptions \eqref{assG}, \eqref{assFGan}. Let $\tilde\xi\in\mathbb{R}^n$. Let assumption \eqref{assxi} hold. Let $\{\varepsilon_j\}_{j \in \mathbb{N}}$ be a sequence in $]0,\epsilon_0[$ converging to $0$. Let $\{u_j\}_{j \in \mathbb{N}}$ be a sequence of functions such that
\begin{align}
& u_j\in  C^{m,\alpha}_{\mathrm{loc}}(\mathrm{cl} \mathbb{S}[\Omega^h_{p,\varepsilon_j}]^-,\mathbb{R}^n)\quad \forall j \in \mathbb{N}\,, \label{eq:uniqueps1}\\
& \text{$u_j$ solves problem \eqref{bvp:nltraceleps} with $\epsilon \equiv \varepsilon_j \quad \forall j \in \mathbb{N}$}\,, \label{eq:uniqueps2}\\
& \text{$\lim_{j \to \infty}u_j(p+\varepsilon_j \,\cdot)_{|\partial\Omega^h}=\tilde{\xi}$ in $C^{m-1,\alpha}(\partial \Omega^h,\mathbb{R}^n)$}\,. \label{eq:uniqueps3} 
\end{align}
Here $u_j(p+\varepsilon_j \,\cdot)_{|\partial\Omega^h}$ denotes the map from $\partial \Omega^h$ to $\mathbb{R}^n$ which takes $t$ to $u_j(p+\varepsilon_j t)$. Then there exists $j_0 \in \mathbb{N}$ such that
\[
u_j= u(\varepsilon_j,\cdot) \qquad  \forall j\in\mathbb{N}\text{ such that } j \geq j_0\,.
\]
\end{thm}
{\bf Proof.} 
Let $\epsilon_1$ be as in Theorem \ref{thm:Lmbd}. By conditions \eqref{eq:uniqueps1}, \eqref{eq:uniqueps2}, and  Proposition \ref{prop:biju}, for each $j \in \mathbb{N}$ there exists a unique pair $(\theta_j,\xi_j)$ in $C^{m-1,\alpha}(\partial \Omega^h,\mathbb{R}^n)_0\times \mathbb{R}^n$ such that
\begin{equation}\label{eq:uniqueps9}
u_j(x)-{Bq^{-1}}x=u_\#[\varepsilon_j,\theta_j,\xi_j](x) \qquad \forall x \in \mathrm{cl}\mathbb{S}[\Omega^h_{p,\varepsilon_j}]^{-}\,. 
\end{equation}
Then to show the validity of the Theorem, it will be enough to prove that 
\begin{equation}\label{eq:uniqueps4}
\lim_{j \to \infty}(\theta_j,\xi_j)=(\tilde{\theta},\tilde{\xi}) \qquad \text{in $C^{m-1,\alpha}(\partial \Omega^h,\mathbb{R}^n)_0\times\mathbb{R}^n$}\,.
\end{equation}
Indeed, if we denote by $\mathcal{U}$ the neighbourhood of Theorem \ref{thm:Lmbd}, the limiting relation in \eqref{eq:uniqueps4} implies that there exists $j_0 \in \mathbb{N}$ such that
$(\varepsilon_j,\theta_j,\xi_j)\in ]0,\epsilon_1[\times \mathcal{U}$
for all $j \geq j_0$ and thus Theorem \ref{thm:Lmbd} would imply that $(\theta_j,\xi_j)=(\Theta[\varepsilon_j],\Xi[\varepsilon_j])$
for all $j \geq j_0$, and that accordingly the Theorem holds (cf.~Definition \ref{def:ueps}). Thus we now turn to the proof of the limit in \eqref{eq:uniqueps4}. We note that equation $\Lambda[\epsilon,\theta,\xi]=0$ can be rewritten in the following form
\begin{eqnarray}
\label{eq:uniqueps5a}
\nonumber
\lefteqn{\frac{1}{2}\theta(t)+w_{\ast}[\omega,\theta](t)+\epsilon^{n-1} \int_{\partial \Omega^h}\sum_{l=1}^n \theta_{l}(s)T(\omega,D R_{n,\omega}^{q,l}(\epsilon(t-s)))\nu_{\Omega^h}(t)\,d\sigma_s}\\
\nonumber
\lefteqn{-D_uG(t,\tilde{\xi})\Bigl(\epsilon v[\omega,\theta](t)+\epsilon^{n-1}\int_{\partial \Omega^h}R_{n,\omega}^q(\epsilon(t-s))\theta(s)\,d\sigma_s+ \xi\Bigr)}\\
\nonumber
&&= G\Bigl(t,\epsilon v[\omega,\theta](t)+\epsilon^{n-1}\int_{\partial \Omega^h}R_{n,\omega}^q(\epsilon(t-s))\theta(s)\,d\sigma_s+\epsilon {Bq^{-1}}t + \xi\Bigr)\\
\nonumber
&&\quad -D_uG(t,\tilde{\xi})\Bigl(\epsilon v[\omega,\theta](t)+\epsilon^{n-1}\int_{\partial \Omega^h}R_{n,\omega}^q(\epsilon(t-s))\theta(s)\,d\sigma_s+ \xi\Bigr)\\
&&\quad -T(\omega,{Bq^{-1}})\nu_{\Omega^h}(t)\qquad\qquad\qquad\qquad\qquad\qquad \qquad \forall t \in \partial \Omega^h\,,
\end{eqnarray}
for all $(\epsilon,\theta,\xi)$ in $]-\epsilon_0,\epsilon_0[\times C^{m-1,\alpha}(\partial \Omega^h,\mathbb{R}^n)_0\times \mathbb{R}^n$. We define the map $N$ from $]-\epsilon_1,\epsilon_1[\times C^{m-1,\alpha}(\partial \Omega^h,\mathbb{R}^n)_0\times \mathbb{R}^n$ to $C^{m-1,\alpha}(\partial \Omega^h,\mathbb{R}^n)$ by setting $N[\epsilon,\theta,\xi]$ equal to the left-hand side of the equality in \eqref{eq:uniqueps5a}. By the proof of Theorem \ref{thm:Lmbd}, $N$ is real analytic. Since $N[\epsilon,\cdot,\cdot]$ is linear for all $\epsilon \in ]-\epsilon_1,\epsilon_1[$, we have
\[
N[\epsilon,\theta,\xi]=\partial_{(\theta,\xi)}N[\epsilon,\tilde{\theta},\tilde{\xi}](\theta,\xi)
\]
for all $(\epsilon,\theta,\xi) \in ]-\epsilon_1,\epsilon_1[\times C^{m-1,\alpha}(\partial \Omega^h,\mathbb{R}^n)_0\times \mathbb{R}^n$, and the map from $]-\epsilon_1,\epsilon_1[$ to $\mathcal{L}(C^{m-1,\alpha}(\partial \Omega^h,\mathbb{R}^n)_0\times \mathbb{R}^n,C^{m-1,\alpha}(\partial \Omega^h,\mathbb{R}^n))$ which takes $\epsilon$ to $N[\epsilon,\cdot,\cdot]$ is real analytic. Moreover,
\[
N[0,\cdot,\cdot]=\partial_{(\theta,\xi)}\Lambda[0,\tilde{\theta},\tilde{\xi}](\cdot,\cdot)\,.
\]
Thus the proof of Theorem \ref{thm:Lmbd} implies that $N[0,\cdot,\cdot]$ is also a linear homeomorphism. As is well known, the set of linear homeomorphisms from $C^{m-1,\alpha}(\partial \Omega^h,\mathbb{R}^n)_0\times \mathbb{R}^n$ to $C^{m-1,\alpha}(\partial \Omega^h,\mathbb{R}^n)$ is open in $\mathcal{L}(C^{m-1,\alpha}(\partial \Omega^h,\mathbb{R}^n)_0\times \mathbb{R}^n,C^{m-1,\alpha}(\partial \Omega^h,\mathbb{R}^n))$ and the map which takes a linear invertible operator to its inverse is real analytic (cf., \textit{e.g.}, Hille and Phillips \cite[Theorems 4.3.2 and 4.3.4]{HiPh57}). Therefore there exists $\epsilon_2 \in]0,\epsilon_1[$ such that the map $\epsilon \mapsto N[\epsilon,\cdot,\cdot]^{(-1)}$ is real analytic from $]-\epsilon_2,\epsilon_2[$ to $\mathcal{L}(C^{m-1,\alpha}(\partial \Omega^h,\mathbb{R}^n), C^{m-1,\alpha}(\partial \Omega^h,\mathbb{R}^n)_0\times\mathbb{R}^n)$. We now denote by $S[\epsilon,\theta,\xi]$ the function defined by the right-hand side of equation \eqref{eq:uniqueps5a}. Then equation $\Lambda[\epsilon,\theta,\xi]=0$ (or equivalently equation \eqref{eq:uniqueps5a}) can be rewritten in the following form,
\begin{equation}\label{eq:uniqueps6}
(\theta,\xi)=N[\epsilon,\cdot,\cdot]^{(-1)}[S[\epsilon,\theta,\xi]]\,,
\end{equation}
for all $(\epsilon,\theta,\xi) \in ]-\epsilon_2,\epsilon_2[\times C^{m-1,\alpha}(\partial \Omega^h,\mathbb{R}^n)_0\times \mathbb{R}^n$. Next we note that the equality in \eqref{eq:uniqueps9} and the definition of $u_\sharp[\cdot,\cdot,\cdot]$ in \eqref{eq:biju2} imply that 
\[
S[\varepsilon_j,\theta_j,\xi_j](t)=G(t,u_j(p+\varepsilon_jt)) -D_uG(t,\tilde\xi)(u_j(p+\varepsilon_jt)-\varepsilon_j Bq^{-1}t)-T(\omega,Bq^{-1})\nu_{\Omega^h}(t)\qquad\forall t\in\partial\Omega^h\,,\ j\in\mathbb{N}.
\]
Then, by condition \eqref{eq:uniqueps3}, and by the real analyticity of $F_G$, and by standard calculus in Banach space we deduce that 
\begin{equation}\label{eq:uniqueps8}
\lim_{j \to \infty}S[\varepsilon_j,\theta_j,\xi_j]=G(\cdot,\tilde\xi)-D_uG(\cdot,\tilde\xi)\tilde\xi-T(\omega,Bq^{-1})\nu_{\Omega^h}=S[0,\tilde{\theta},\tilde{\xi}] 
\end{equation} 
in $C^{m-1,\alpha}(\partial \Omega^h,\mathbb{R}^n)$.
Then by equality \eqref{eq:uniqueps6}, and by the limit in \eqref{eq:uniqueps8}, and by the real analyticity of the map which takes $\epsilon$ to $N[\epsilon,\cdot,\cdot]^{(-1)}$, and by the bilinearity and continuity of the operator from $\mathcal{L}(C^{m-1,\alpha}(\partial \Omega^h,\mathbb{R}^n),C^{m-1,\alpha}(\partial \Omega^h,\mathbb{R}^n)_0\times \mathbb{R}^n)\times C^{m-1,\alpha}(\partial \Omega^h,\mathbb{R}^n)$ to $C^{m-1,\alpha}(\partial \Omega^h,\mathbb{R}^n)_0\times \mathbb{R}^n$, which takes a pair $(T_1,T_2)$ to $T_1[T_2]$, we conclude that the limit in \eqref{eq:uniqueps4} holds. Thus the proof is complete.
 \hfill $\Box$
\vspace{\baselineskip}

\section*{Acknowledgements}

The research  of M.~Dalla Riva was supported by {\it FEDER} funds through {\it COMPETE}--Operational Programme Factors of Competitiveness (``Programa Operacional Factores de Competitividade'') and by Portuguese funds through the {\it Center for Research and Development in Mathematics and Applications} (University of Aveiro) and the Portuguese Foundation for Science and Technology (``FCT--Funda{\c c}{\~a}o para a Ci\^encia e a Tecnologia''), within project PEst-C/MAT/UI4106/2011 with the COMPETE  number FCOMP-01-0124-FEDER-022690. The research was also supported  by the Portuguese Foundation for Science and Technology (``FCT--Funda{\c c}{\~a}o para a Ci\^encia e a Tecnologia'')  with the research grant SFRH/BPD/64437/2009. The research of P.~Musolino was supported by the ``Accademia Nazionale dei Lincei'' through a scholarship ``Royal Society''. Part of the work was done while P.~Musolino was visiting the Centro de Investiga\c{c}\~{a}o e Desenvolvimento em Matem\'{a}tica e Aplica\c{c}o\~{e}s of the Universidade de Aveiro. P.~Musolino wishes to thank the Centro de Investiga\c{c}\~{a}o e Desenvolvimento em Matem\'{a}tica e Aplica\c{c}o\~{e}s, and in particular Prof.~L.~P.~Castro and Dr.~M.~Dalla Riva, for the kind hospitality.

\end{document}